\documentclass[11pt]{amsart}

\usepackage{amssymb,amsthm,amsmath,hyperref,amsrefs,verbatim}

\usepackage[margin=1.25in]{geometry}

\numberwithin{equation}{section}

\input{xy}
\xyoption{all}

\newtheorem{theorem}{Theorem}[subsection]
\newtheorem{corollary}[theorem]{Corollary}
\newtheorem{lemma}[theorem]{Lemma}
\newtheorem{proposition}[theorem]{Proposition}

\theoremstyle{remark}
\newtheorem{remark}[theorem]{Remark}

\newcommand{\N}{\mathbb{N}}
\newcommand{\Q}{\mathbb{Q}}
\newcommand{\Z}{\mathbb{Z}}
\newcommand{\R}{\mathbb{R}}
\newcommand{\C}{\mathbb{C}}
\newcommand{\CCu}{\mathbf{Cu}}
\newcommand{\Cu}{\mathrm{Cu}}

\newcommand{\W}{\mathrm{W}}
\newcommand{\T}{\mathrm{T}}

\newcommand{\cuntzle}{\preccurlyeq}
\renewcommand{\epsilon}{\varepsilon}
\renewcommand{\leq}{\leqslant}
\renewcommand{\geq}{\geqslant}
\newcommand{\rank}{\mathrm{rank}\ }

\title{Classification of inductive limits of 1-dimensional NCCW complexes}
\author{Leonel Robert}

\address{Department of Mathematics, University of Louisiana at Lafayette, Lafayette USA.}

\email{lrobert@louisiana.edu}

\begin{document}
\begin{abstract}
A classification result is obtained for the C$^*$-algebras that are
(stably isomorphic to) inductive limits of 1-dimensional noncommutative CW complexes with trivial $\mathrm K_1$-group.
The classifying functor $\Cu^\sim$ is defined in terms of the Cuntz semigroup of the unitization of the algebra.
For the simple C$^*$-algebras covered by the classification, $\Cu^\sim$
reduces to the ordered $\mathrm K_0$-group, the cone of traces, and the pairing between them. As an application of the classification, it is shown that the crossed
products by a quasi-free action $\mathcal O_2\rtimes_\lambda \R$ are all isomorphic for a dense set of positive irrational numbers $\lambda$.
\end{abstract}

\maketitle

\setcounter{tocdepth}{1}
\tableofcontents

\section{Introduction}
This  paper is a contribution to the classification of C$^*$-algebras by means of
data invariant under approximately inner automorphisms. The classifying functor used here, denoted by $\Cu^\sim$,  
is a variation on the Cuntz semigroup functor that looks at the unitizations of the C$^*$-algebras.
The C$^*$-algebras being classified are those expressible -- up to stable isomorphism -- as inductive limits of 1-dimensional noncommutative
CW complexes with trivial $\mathrm K_1$-group. 
The following theorem is the main result of the paper; see Sections \ref{prelims} and \ref{Cusim} for the relevant definitions.

\begin{theorem}\label{main}
Let $A$ be either a 1-dimensional noncommutative \emph{CW} complex with trivial $\mathrm K_1$-group, or a sequential inductive limit of such C$^*$-algebras,
or a C$^*$-algebra stably isomorphic to one such inductive limit. 
Let $B$ be a stable rank one C$^*$-algebra. Then for every morphism in the category $\CCu$
\[
\alpha\colon \Cu^\sim(A)\to \Cu^\sim(B)
\]
such that $\alpha([s_A])\leq [s_B]$, where $s_A\in A^+$ and $s_B\in B^+$ are strictly
positive elements, there exists a homomorphism 
\[
\phi\colon A\to B 
\]
such that $\Cu^\sim(\phi)=\alpha$.
Moreover, $\phi$ is unique up to approximate unitary equivalence.
\end{theorem}

The above  classification result continues, and extends, the work of  Ciuperca and Elliott in \cite{ciuperca-elliott},
and  of  Ciuperca, Elliott, and Santiago in \cite{ciuperca-elliott-santiago}. In \cite{ciuperca-elliott} and \cite{ciuperca-elliott-santiago} the classifying functor
is the Cuntz semigroup, while the range of the C$^*$-algebra
$A$ is the class of approximate interval C$^*$-algebras in \cite{ciuperca-elliott}, and  of approximate tree C$^*$-algebras
in \cite{ciuperca-elliott-santiago}. 
These classification results have the following distinctive  characteristics:
\begin{enumerate}
\item
the main theorem classifies homomorphisms from a certain class of domain algebras into
arbitrary stable rank one C$^*$-algebras (in \cite[Theorem 2]{robert-santiago},  it is shown that
the class of codomain algebras can  be allowed to be somewhat larger than the stable rank one C$^*$-algebras);

\item
in contrast to the generality of the codomain algebra, the domain algebra belongs to a relatively special class (up
to stable isomorphism): 
approximate interval C$^*$-algebras, approximate tree C$^*$-algebras, or, as in this paper, 
inductive limits of 1-dimensional noncommutative CW (NCCW) complexes with trivial $\mathrm K_1$-group;

\item
the isomorphism theorem (i.e., the classification
of C$^*$-algebras) is obtained from the classification of homomorphisms  by a well-known approximate intertwining argument; 
in this case both C$^*$-algebras being compared must belong to the class of domain algebras covered
by the homomorphism theorem (see Corollary \ref{isomorphism});

\item
a certain uniform continuity of the classification is also obtained: if two homomorphisms are close at the level of the invariant, then they
are close -- up to approximate unitary equivalence, on finite sets -- as C$^*$-algebra homomorphisms (see Theorem \ref{uniform});

\item
neither the domain nor the codomain  algebra is required to be simple.
\end{enumerate}

For simple C$^*$-algebras, the connection with the standard Elliott invariant (i.e., K-theory and traces) is made as follows: 
If a  C$^*$-algebra $A$ is simple and belongs to the class of domain algebras covered by Theorem \ref{main}, 
then the invariant $\Cu^\sim(A)$ can be computed in terms of the ordered group $\mathrm K_0(A)$, the cone of traces of $A$, and the pairing between them. 
It then follows by our classification  that the simple C$^*$-algebras that are (stably isomorphic to) inductive limits of 1-dimensional NCCW complexes with trivial $\mathrm K_1$-group are classified by their ordered $\mathrm K_0$-group, their cone of traces, and the pairing between them (a suitable scale is also
needed; see Corollary \ref{classsimple}).

Among the simple C$^*$-algebras included in our classification are
\begin{enumerate}
\item
the Jiang-Su algebra $\mathcal Z$ (a finite, nuclear, C$^*$-algebra not of type I, with the same Elliott invariant as $\C$);

\item
the simple projectionless C$^*$-algebras classified by Razak in \cite{razak};

\item
the simple C$^*$-algebras with $\mathrm K_1=0$ described by Elliott in \cite[Theorem 2.2]{elliott2};
their Elliott invariants exhaust all the pairs $(G,C)$, with $G$ a torsion free countable ordered abelian group
and  $C$ a non-zero topological cone with a compact base that is a metrizable Choquet simplex, where
a pairing $G\times C\to \R$ is given that is weakly unperforated (see \cite{elliott2} for details);

\item
the simple inductive limits of splitting interval algebras classified by Jiang and Su in \cite{jiang-su} (including the ordinary interval algebra
case considered by Elliott in \cite{elliott92}).  The non-simple inductive limits of splitting interval algebras classified by Santiago in \cite{santiago} are also included in our classification.
\end{enumerate}
As an application of the classification, we show that the Jiang-Su algebra embeds unitally into any non-elementary unital simple exact C$^*$-algebra of stable rank one, with a unique tracial state,
and with strict comparison of positive elements. Furthermore, the embedding is unique up to approximate unitary equivalence.
This property characterizes the Jiang-Su algebra up to isomorphism (see Proposition \ref{Zembeddings}). 

In \cite{dean}, Dean studied the  crossed products $\mathcal O_2\rtimes_\lambda \R$ obtained from a  quasi-free action of $\R$ on the Cuntz algebra $\mathcal O_2$. He showed that
for $\lambda$ in a  dense $G_\delta$ subset of $\R^+$ containing the rational
numbers these crossed products are inductive limits of 1-dimensional NCCW complexes. A careful examination of his result reveals that the building blocks of the inductive limits have trivial $\mathrm K_1$-group. (Note that this says more than just that the $\mathrm K_1$-group of the limit is zero.)
Thus, such crossed products are included in the classification given here. As an application of the classification, we obtain that for $\lambda$ in a dense 
subset of $\R^+\backslash \Q$ of second Baire category the crossed products
$\mathcal O_2\rtimes_\lambda \R$
are all isomorphic to each other.  In fact, they are isomorphic to the  unique (up to isomorphism) simple, stable, projectionless C$^*$-algebra with a unique trace  
expressible as an inductive limit of 1-dimensional NCCW complexes with trivial $\mathrm K_1$-group. 
In \cite{jacelon}, Jacelon 
 sets forth the study of this C$^*$-algebra
as a stably finite analogue of the Cuntz algebra $\mathcal O_2$ and as a model for what a non-unital
strongly self-absorbing C*-algebra should be. This C$^*$-algebra is likely to play
a significant role in the classification of projectionless simple nuclear C$^*$-algebras, akin to the role played by the unital strongly self-absorbing 
C$^*$-algebras.

Let us say a few words  about the proof of Theorem \ref{main}.  It is shown in  \cite{ciuperca-elliott-santiago}  that the classification
of homomorphisms into a stable rank one C$^*$-algebra by means of the functor $\Cu$
 has certain permanence properties with respect to the domain algebra. Namely, if the classification is possible for a given collection of domain algebras, it is also possible  for their finite direct sums, sequential inductive limits, and for algebras stably isomorphic to the ones in the given collection. 
In Theorem \ref{permanence} below the same permanence results are obtained  for $\Cu^\sim$. Furthermore, the operations of adding and removing a unit are added to the list of transformations of the domain algebra. Theorem \ref{main} is then proved by showing that all
1-dimensional NCCW complexes with trivial $\mathrm K_1$-group may  be gotten by starting with the algebra $\mathrm C_0(0,1]$ and combining the operations of direct sum, adding or removing a unit, and passing to  stably isomorphic algebras. This reduces proving Theorem \ref{main} to the special case $A=\mathrm C_0(0,1]$. In this case
the theorem is essentially a corollary of Ciuperca and Elliott's classification \cite{ciuperca-elliott}.

The scope of the methods used in this paper is currently limited to domain algebras that
have trivial $\mathrm K_1$-group. This is so because  the functors $\Cu$ and $\Cu^\sim$ fail in general to account for the $\mathrm K_1$-group of the algebra. It seems plausible that a suitable enlargement of these functors with $\mathrm K_1$-type data may give a classification of all 1-dimensional NCCW complexes and their inductive limits.

This paper is organized as follows: Section \ref{prelims} contains preliminaries about the Cuntz semigroup and NCCW complexes; in Section \ref{Cusim}
the functor $\Cu^\sim$ is introduced and studied; Section \ref{specialcases}
contains  proofs of various special cases of Theorem \ref{main} and serves as a warm-up for the proof of the general case; Section \ref{proofofmain}
is dedicated to the proof of Theorem \ref{main}; in Section \ref{simple} a computation is given of $\Cu^\sim(A)$ for $A$ simple
and belonging to the class classified here; the classification of simple C$^*$-algebras in terms of the ordered $\mathrm K_0$-group and the cone of traces is derived from this computation;
in Section \ref{crossed} the crossed products $\mathcal O_2\rtimes_\lambda \R$ are discussed. 
    
\textbf{Acknowledgments.} I am grateful to Alin Ciuperca, George A. Elliott, and Luis Santiago for fruitful discussions on the 
role of the Cuntz semigroup in the classification of C$^*$-algebras, to Luis Santiago for pointing out a gap in the original reduction algorithm of Section \ref{proofofmain}, and to Mikael R\o rdam for sharing with me his Proposition \ref{proprordam}.

\section{Preliminary definitions and results}\label{prelims}
\subsection{The functor $\Cu$}\label{prelimsCu}
Let us briefly review the definition of the Cuntz semigroup of a C$^*$-algebra and of the
functor $\Cu$ (see \cite{ara-perera-toms} for a more detailed exposition). Here we use the positive elements picture of $\Cu$. 
In \cite{coward-elliott-ivanescu}, an alternative approach to $\Cu$ is given  which makes use of Hilbert C$^*$-modules over the C$^*$-algebra
rather than positive elements, but we will not rely on it here.

Let $A$ be a C$^*$-algebra. Let $A^+$ denote the positive elements of $A$. Given $a,b\in A^+$ we say that $a$ is Cuntz smaller than $b$, and denote this by  $a\cuntzle b$,
if $d_nbd_n^*\to a$ for some sequence $(d_n)$ in $A$. We say that $a$ is Cuntz equivalent to $b$ if  $a\cuntzle b$ and $b\cuntzle a$; in 
this case we write $a\sim b$.

Let $\Cu(A)$ denote the set $(A\otimes \mathcal K)^+/\!\!\sim$ of Cuntz equivalence classes of positive elements of $A\otimes\mathcal K$.
For $a\in (A\otimes\mathcal K)^+$, let us denote the Cuntz class of $a$ by $[a]$.  The preorder $\cuntzle$ defines an order on $\Cu(A)$:
\[
[a]\leq [b]\hbox { if }a\cuntzle b.
\]
We also endow $\Cu(A)$ with an addition operation by setting
\[
[a]+[b]:=[a'+b'],
\]
where
$a',b'\in (A\otimes\mathcal K)^+$ are orthogonal to each other and Cuntz equivalent to $a$ and $b$ respectively (the choices of $a'$ and $b'$
do not affect the Cuntz class of their sum). We regard $\Cu(A)$ as an ordered semigroup and we call it
the Cuntz semigroup of $A$.

Let $\phi\colon A\to B$ be a homomorphism between C$^*$-algebras. 
Let us continue to denote by $\phi$ the homomorphism $\phi\otimes \mathrm{id}$ from $A\otimes\mathcal K$ to $B\otimes\mathcal K$. (This convention
will apply throughout the paper.) The homomorphism $\phi$
 induces an ordered semigroup map $\Cu(\phi)\colon \Cu(A)\to \Cu(B)$ given by
\[
\Cu(\phi)([a]):=[\phi(a)].
\]

In \cite[Theorem 1]{coward-elliott-ivanescu}, Coward, Elliott and Ivanescu show that $\Cu(\cdot)$
is a functor from the category of C$^*$-algebras to a certain sub-category of the category of ordered abelian semigroups.
Let us recall the definition of this sub-category,  which we shall  denote by $\CCu$.

The definition of the category $\CCu$ relies on the compact containment relation.
For $x$ and $y$ elements of an ordered set,  we say that $x$ is compactly contained in $y$, and denote this by $x\ll y$,
if for every increasing sequence $(y_n)$ with $y\leq \sup_n y_n$, there exists an index $n_0$ such that $x\leq y_{n_0}$. A sequence
$(x_n)$ is called rapidly increasing if $x_n\ll x_{n+1}$ for all $n$.
An ordered abelian semigroup with zero $S$ is an object of $\CCu$ if
\begin{enumerate}
\item [O1]
increasing sequences in $S$ have suprema,

\item [O2]
for every $x\in S$ there exists a rapidly increasing sequence $(x_n)$ such $x=\sup_n x_n$,

\item [O3]
if $x_i\ll y_i$, $i=1,2$, then $x_1+x_2\ll y_1+y_2$,

\item [O4]
$\sup_n (x_n+y_n)=\sup_n(x_n)+\sup_n(y_n)$, for $(x_n)$ and $(y_n)$ increasing sequences.
\end{enumerate}

A map $\alpha\colon S\to T$ is a morphism of $\CCu$ if
\begin{enumerate}
\item[M1] $\alpha$ is additive and maps 0 to 0,

\item[M2] $\alpha$ is order preserving,

\item [M3]
$\alpha$ preserves the  suprema of increasing sequences,

\item [M4]
$\alpha$ preserves the relation of compact containment.
\end{enumerate}

It is shown in \cite[Theorem 2]{coward-elliott-ivanescu} that the category $\CCu$ is closed under sequential inductive limits, 
and that the functor $\Cu$ preserves sequential inductive limits. For later use, we will need the following characterization of inductive limits 
in the category $\CCu$:
Given an inductive system $(S_i,\alpha_{i,j})_{i,j\in \N}$ in the category $\CCu$, an object $S$ in $\CCu$, and morphisms $\alpha_{i,\infty}\colon S_i\to S$ 
such that $\alpha_{j,\infty}\circ \alpha_{i,j}=\alpha_{i,\infty}$ for all $i\leq j$, the object $S$ is the inductive limit of the $S_i$s if
\begin{enumerate}
\item [L1]
the set $\bigcup_i \alpha_{i,\infty}(S_i)$ is dense in $S$ (in the sense that every element of $S$ is the supremum
of a rapidly increasing  sequence of elements in $\bigcup_i \alpha_{i,\infty}(S_i)$),

\item [L2]
 for each $x,y\in S_i$ such that $\alpha_{i,\infty}(x)\leq \alpha_{i,\infty}(y)$, and $x'\ll x$,
there exists $j$ such that $\alpha_{i,j}(x')\leq \alpha_{i,j}(y)$.
\end{enumerate}
(It is easy to see that these properties hold for the inductive limit as constructed in \cite{coward-elliott-ivanescu}, and that
they imply in turn the universal property of an inductive limit in the category $\CCu$.)

The order on $\Cu(A)$ is part of its structure and in general is not determined by the addition operation.
In one special situation, however, the order coincides with the algebraic order. We say that an element $e\in \Cu(A)$
is compact if $e\ll e$. The following lemma is a  simple consequence of \cite[Lemma 7.1 (i)]{rordam-winter}.

\begin{lemma}\label{algebraicorder}
Let  $e\in \Cu(A)$ be compact and suppose that $e\leq [a]$ for some $[a]\in \Cu(A)$. Then there exists $[a']\in \Cu(A)$
such that $e+[a']=[a]$. 
\end{lemma}

We will sometimes make use of the ordered sub-semigroup $\W(A)$ of $\Cu(A)$. This ordered 
semigroup is  composed of 
the Cuntz equivalence classes $[a]\in \Cu(A)$ such that $a\in \bigcup_{n=1}^\infty \mathrm M_n(A)$.
It is known that $\W(A)$ is dense in $\Cu(A)$, i.e.,
every element of $\Cu(A)$ is the supremum of a rapidly increasing sequence of elements in $\W(A)$. Indeed, if $[a]\in \Cu(A)$
then -- by the proof of \cite[Theorem 1]{coward-elliott-ivanescu} -- $[(a-\frac{1}{n})_+]$, with $n=1,2,\dots$, is a rapidly increasing sequence in $\W(A)$ with supremum $[a]$.

In the following sections we will make use of some properties of the Cuntz semigroup that hold specifically for stable rank one C$^*$-algebras. 
Let us recall them here.

The following proposition is due to Coward, Elliott, and Ivanescu
(\cite[Theorem 3]{coward-elliott-ivanescu}; see \cite[Proposition 1]{ciuperca-elliott-santiago} for the present formulation).
\begin{proposition} (\cite[Proposition 1]{ciuperca-elliott-santiago})\label{sr1iso}
Let $A$ be a C$^*$-algebra of stable rank one and $a,b\in A^+$. Then $a\cuntzle b$
if and only if there is $x\in A$ such that $a=x^*x$ and $xx^*\in \overline{bAb}$.
\end{proposition}

Let us say that the ordered (abelian) semigroup $S$ has weak cancellation if $x+z\ll y+z$
implies $x\leq y$, for $x,y,z\in S$. If $S=\Cu(A)$ and $[p]$ is the Cuntz class of a projection,
then $[p]$ is a compact element of $\Cu(A)$, i.e., $[p]\ll [p]$ (as a consequence of the fact proved in \cite[Theorem 1]{coward-elliott-ivanescu} (page 168) that the abstract and concrete notions of compact containment coincide). It can be shown from this that
if $\Cu(A)$ has weak cancellation then $[a]+[p]\leq [b]+[p]$ implies $[a]\leq [b]$. In this
case we say that $\Cu(A)$ has cancellation of projections.
The following proposition is due to R\o rdam and Winter (\cite{rordam-winter}). A slightly weaker statement is  due to Elliott (\cite{elliott}).
\begin{proposition} (\cite[Proposition 4.2, Theorem 4.3]{rordam-winter}.)\label{weak-cancellation}
Let $A$ be a C$^*$-algebra of stable rank one. Then $\Cu(A)$ has weak cancellation.
In particular, $\Cu(A)$ has cancellation of projections.
\end{proposition}

\subsection{Noncommutative CW complexes.}\label{NCCWprelims}
Following \cite{elp}, let us say that the C$^*$-algebra $A$ is a 1-dimensional noncommutative CW complex (NCCW complex) if it is the pull-back C$^*$-algebra in a diagram of the form
\begin{align}\label{nccw}
\xymatrix{
A\ar[r]^-{\pi_1}\ar[d]^{\pi_2} & \mathrm C([0,1],F)\ar[d]^{\mathrm{ev}_0\oplus \mathrm{ev}_1}\\
E\ar[r]^-{\phi} & F\oplus F.
}
\end{align}
Here $E$ and $F$ are finite dimensional C$^*$-algebras and
$\mathrm{ev}_0$ and $\mathrm{ev}_1$ are the evaluation maps at 0 and 1.
As in \cite{dean} -- and unlike \cite{elp} -- we do not assume that the homomorphism 
$\phi$ is unital. (Of course, we are using the expression``CW complex" to mean ``finite CW complex".)

Let us introduce some notation for the data defining $A$.
We may write
\begin{align*}
E &= \mathrm M_{e_1}(\C)\oplus \mathrm M_{e_2}(\C)\oplus \dots \oplus\mathrm M_{e_k}(\C),\\
F &= \mathrm M_{f_1}(\C)\oplus \mathrm M_{f_2}(\C)\oplus \dots \oplus\mathrm M_{f_l}(\C).
\end{align*}
We have $\phi=(\phi_0,\phi_1)$, with $\phi_0,\phi_1\colon E\to F$. 
We may identify $\mathrm K_0(E)$ with $\Z^k$ and $\mathrm K_0(F)$ with $\Z^l$.
Let us  denote by  $Z^{\phi_0}$ and $Z^{\phi_1}$ the $l\times k$ integer matrices associated to the maps
\[
\mathrm K_0(\phi_0)\colon \Z^{k}\to \Z^{l}\hbox{ and }\mathrm K_0(\phi_1)\colon \Z^{k}\to \Z^{l}.
\]

\begin{remark}\label{stablyNCCW} The following observations will be useful.

(i) The vectors $(e_j)_{j=1}^k$ and $(f_i)_{i=1}^l$ determine the algebras $E$ and $F$
uniquely up to isomorphism. These vectors  together with the matrices
$Z^{\phi_0}$ and $Z^{\phi_1}$ determine the maps $\phi_0$ and $\phi_1$ up to unitary equivalence. In consequence, the data $(e_j)_{j=1}^k$, $(f_i)_{i=1}^l$, $Z^{\phi_0}$ and $Z^{\phi_1}$ determine $A$  up to isomorphism.

(ii) Given the data  $(e_j)_{j=1}^k$, $(f_i)_{i=1}^l$, $Z^{\phi_0}$ and $Z^{\phi_1}$,
there is a 1-dimensional NCCW complex that attains them if and only if
\[
Z^{\phi_0}(e_j)\leq (f_i) \hbox{ and }Z^{\phi_1}(e_j)\leq (f_i).
\]

(iii) The operation of tensoring by $\mathcal K$ preserves C$^*$-algebra pull-backs. Thus, tensoring by $\mathcal K$ in the diagram \eqref{nccw} we can deduce that the matrices
$Z^{\phi_0}$ and $Z^{\phi_1}$ alone determine $A$  up to stable isomorphism (since these matrices uniquely determine the stabilized homomorphisms $\phi_0\otimes\mathrm{id},\phi_1\otimes\mathrm{id}\colon E\otimes\mathcal K\to F\otimes\mathcal K$  up to unitary equivalence).
\end{remark}

\subsection{Approximate unitary equivalence}
Given two homomorphisms $\phi\colon A\to B$ and $\psi\colon A\to B$ let us say that they are approximately
unitarily equivalent if there exists a net $(u_\lambda)_{\lambda\in \Lambda}$ of unitaries in $B^\sim$ such that
$u_\lambda\phi(x)u_\lambda^*\to \psi(x)$ for all $x\in A$.  Notice that without loss of generality, the index set $\Lambda$ can be chosen to be the pairs $(F,\epsilon)$, with $F$ finite subset of $A$, $\epsilon>0$, and with the order $(F,\epsilon)\precsim (F',\epsilon')$ if $F\subseteq F'$ and $\epsilon'\leq \epsilon$.

\begin{proposition}\label{sr1hered}
(i) Let $\phi,\psi\colon A\to B$ be approximately unitarily equivalent homomorphisms, with
$B$ of stable rank one. Suppose that $\phi(A),\psi(A)\subseteq C$, with $C$ a closed hereditary subalgebra
of $B$. Then $\phi$ and $\psi$ are approximately unitarily equivalent as homomorphisms from $A$
to $C$.

(ii) Let $C$ be a stable rank one C$^*$-algebra. Let $\phi,\psi\colon A\oplus B\to C$ be homomorphisms such 
that $\phi|_A$ is approximately unitarily equivalent to $\psi|_A$ and $\phi|_B$ is approximately unitarily
equivalent to $\psi|_B$. Then $\phi$ is approximately unitarily equivalent to $\psi$.
\end{proposition}
\begin{proof}
(i) (Cf. proof of \cite[Lemma 6]{ciuperca-elliott-santiago}.) Let $(u_\lambda)_{\lambda\in\Lambda}$ be a net of unitaries such that $u_\lambda\phi(x)u_\lambda^*\to \psi(x)$ for $x\in A$.
Without loss of generality, we may assume that $\Lambda$ is the set of pairs $(F,\epsilon)$, with $F$ a finite subset of $A$ and $\epsilon>0$.
It suffices to assume that $C$ is the closed hereditary subalgebra generated
by the images of $\phi$ and $\psi$. For this subalgebra, we may choose an approximate unit
$(e_\lambda)_{\lambda\in \Lambda}$ indexed by $\Lambda$ (e.g., for $\lambda=(F,\epsilon)$ choose $e_\lambda$ such that
$\|e_\lambda\phi(x)-\phi(x)\|<\epsilon$ and $\|e_\lambda\psi(x)-\psi(x)\|<\epsilon$ for all $x\in F$).
Set
\[
e_\lambda u_\lambda e_\lambda=z_\lambda.
\]
Then
\begin{align}\label{zlambda1}
z_\lambda\phi(x)z_\lambda^*\to \psi(x)\hbox{ for all }x\in A,\\
|z_\lambda|\phi(x)\to \phi(x)\hbox{ for all }x\in A.\label{zlambda2}
\end{align}
Since $C$ is of stable rank one, for each $\lambda$ there exists a unitary $w_\lambda\in C^\sim$ such that
$\|z_\lambda-w_\lambda|z_\lambda|\|<\epsilon$, where $\lambda=(F,\epsilon)$. 
(The existence of $w_\lambda$ follows immediately from the
definition of stable rank one by considering the polar decomposition of an invertible element sufficiently close
to $z_\lambda$.)
 Then \eqref{zlambda1} and \eqref{zlambda2} imply that
$w_\lambda\phi(x)w_\lambda^*\to \psi(x)$ for all $x\in A$. That is, $\phi$ and $\psi$ are approximately
unitarily equivalent as homomorphisms with codomain $C$.

(ii) (Cf., proof of \cite[Proposition 5 (iii)]{ciuperca-elliott-santiago}.) Let $(u_\lambda)$ and $(v_\kappa)$ be nets of unitaries in $C^\sim$ such that
$u_\lambda\phi(x)u_\lambda^*\to \psi(x)$ for all $x\in A$ and $v_\kappa\phi(x)v_\kappa^*\to \psi(x)$
for all $x\in B$. Choose approximate units $(e_\lambda^\phi)$ and $(e_\lambda^\psi)$ for the hereditary
subalgebras generated by $\phi(A)$ and $\psi(A)$ respectively. Similarly, choose  approximate units $(f_\kappa^\phi)$ and $(f_\kappa^\psi)$
for the hereditary subalgebras generated by $\phi(B)$ and $\psi(B)$. Set
\[
e_\lambda^\phi u_\lambda e_\lambda^\psi+f_\kappa^\psi v_\kappa f_\kappa^{\phi}=z_{\lambda,\kappa}.
\]
Then
\begin{align}\label{zlambdakappa1}
z_{\lambda,\kappa}\phi(x)z_{\lambda,\kappa}^*\to \psi(x)\hbox{ for all }x\in A\oplus B,\\
|z_{\lambda,\kappa}|\phi(x)\to \phi(x)\hbox{ for all }x\in A\oplus B.\label{zlambdakappa2}
\end{align}
Since $C$ is of stable rank one, for each $\lambda$ and $\kappa$ there exists a  unitary $w_{\lambda,\kappa}\in C^\sim$
such that $\|z_{\lambda,\kappa}-w_{\lambda,\kappa}|z_{\lambda,\kappa}|\|<\epsilon,\epsilon'$,
where $\lambda=(F,\epsilon)$ and $\kappa=(F',\epsilon')$. Then \eqref{zlambdakappa1} and \eqref{zlambdakappa2}
imply that $w_{\lambda,\kappa}\phi(x)w_{\lambda,\kappa}^*\to \psi(x)$
for all $x\in A\oplus B$.
\end{proof}

\section{The functor $\Cu^\sim$}\label{Cusim}
\subsection{Definition and properties of $\Cu^\sim$}
Here we define the functor $\Cu^\sim$ and show that it is well behaved with respect to
sequential inductive limits, stable isomorphism, and direct sums,  as long as the C$^*$-algebras
are assumed to be of stable rank one. This is sufficient generality for the applications to classification that we will consider later on, but it raises the question of whether the same properties of $\Cu^\sim$ hold for a larger class of C$^*$-algebras.

Let $A$ be a C$^*$-algebra and let $A^\sim$ denote its unitization.
The definition of  the  ordered semigroup $\Cu^\sim(A)$ is analogous to the definition of 
$\mathrm K_0(A)$ in terms of the Murray-von Neumann semigroup of $A^\sim$ (N.B.: $\Cu^\sim(A)$ is not the enveloping group of the Cuntz semigroup of $A^\sim$.) 
Just as is sometimes done when defining $\mathrm K_0(A)$, we first define
$\Cu^\sim(A)$ for $A$ unital  and then extend the definition to  an arbitrary $A$ by requiring that the sequence
\[
\Cu^\sim(A)\to \Cu^\sim(A^\sim)\to \Cu^\sim(\C)
\]
be exact in the middle. (This raises the question whether $\Cu^\sim$ is half-exact in general. In the stable rank one setting,
this is shown in Proposition \ref{splitexact} below.)

Let $A$ be a unital C$^*$-algebra. Let us define  $\Cu^\sim(A)$ as the ordered semigroup of formal differences $[a]-n[1]$,
with $[a]\in \Cu(A)$ and  $n\in \N$. That is, $\Cu^\sim(A)$ is the quotient of the   semigroup of pairs $([a],n)$, with $[a]\in \Cu(A)$ and $n\in \N$, by  the equivalence relation
$([a],n)\sim ([b],m)$ if
\[
[a]+m[1]+k[1]=[b]+n[1]+k[1],
\]
for some $k\in \N$. The image of $([a],n)$ in this quotient
will be denoted by $[a]-n[1]$.  One obtains an order relation on $\Cu^\sim(A)$ by saying that 
$[a]-n[1]\leq [b]-m[1]$ if for some $k$ the inequality $[a]+m[1]+k[1]\leq [b]+n[1]+k[1]$
holds in $\Cu(A)$.

An alternative picture of $\Cu^\sim(A)$ for $A$ unital which we will find convenient to have is as the ordered semigroup of formal differences $[a]-e$, with $[a]\in \Cu(A)$ and $e\in \Cu(A)$ a compact element.
Since for every compact element $e$ there exists $e'\in \Cu(A)$ such that
$e+e'=n[1]$ for some $n$ (by Lemma \ref{algebraicorder}), this ordered semigroup coincides with the above defined semigroup of formal differences $[a]-n[1]$.
Since the compact elements of $\Cu(A)$ are intrinsically determined by its order structure, $\Cu^\sim(A)$ is completely determined by $\Cu(A)$ if $A$ is unital. A similar argument shows that $\Cu^\sim(A)$ may also be viewed as the 
ordered semigroup of formal differences $[a]-[p]$, with $[a]\in \Cu(A)$ and $p\in A\otimes\mathcal K$ a projection
(because projections are also complemented, i.e., $[p]+[p']=n[1]$ for some projection $p'$ and $n\in \N$).
%

Let us now define the ordered semigroup $\Cu^\sim(A)$ for an arbitrary  C$^*$-algebra $A$. Let $\pi\colon A^\sim\to \C$ denote the quotient map
from the unitization of $A$ onto $\C$. This map induces
a morphism
\[
\Cu(\pi)\colon \Cu(A^\sim)\to \Cu(\C)\cong  \{0,1,\dots,\infty\},
\]
and also a morphism 
\[
\Cu^\sim(\pi)\colon \Cu^\sim(A^\sim)\to \Cu^\sim(\C)\cong  \Z\cup\{\infty\}.
\]
Let us define $\Cu^\sim(A)$ as the subsemigroup of $\Cu^\sim(A^\sim)$ consisting of the elements 
$[a]-n[1]$, with $[a]$ in $\Cu(A^\sim)$ such that $\Cu(\pi)([a])=n<\infty$. With respect to the relative
order structure, $\Cu^\sim(A)$ is an ordered semigroup.

\begin{remark}\label{CusimArem}
The following facts are readily verified.

(i) If $A$ is unital, but we ignore this fact and find $\Cu^\sim(A)$
by unitizing it, we obtain an ordered semigroup (canonically) isomorphic to $\Cu^\sim(A)$ as
defined above.

(ii) As noted above, if $A$ is unital  then $\Cu^\sim(A)$ is completely determined by $\Cu(A)$. This, however,  may not be the case  if $A$  is non-unital (see Remark \ref{Cufails}).
\end{remark}

There is a canonical map from $\Cu(A)$ to the positive elements of $\Cu^\sim(A)$ 
(i.e., the elements greater than or equal to 0 in $\Cu^\sim(A^\sim)$)
given by $[a]\mapsto [a]-0\cdot [1]$.
\begin{lemma} \label{sur-inj}
The ordered semigroup
$\Cu(A)$ is mapped surjectively onto the positive elements of $\Cu^\sim(A)$,
and injectively if $\Cu(A^\sim)$ has cancellation of projections.
\end{lemma}
Note that if $A$ has stable rank one then $\Cu(A^\sim)$ has cancellation of projections
by Proposition \ref{weak-cancellation}.
\begin{proof}
Let $[a]\in \Cu(A^\sim)$ and $n\in \N$ be such that $[\pi(a)]=n$ (we continue to denote by $\pi$ the extension of the quotient map $\pi\colon A^\sim\to \C$
to $A^\sim\otimes\mathcal K$). If $[a]-n[1]\geq 0$ in $\Cu^\sim(A)$, then $[a]+k[1]\geq (n+k)[1]$ for some $k$, the latter inequality taken in
$\Cu(A^\sim)$. It follows by Lemma \ref{algebraicorder} that $[a]+k[1]=(n+k)[1]+[a']$ for some $[a']$ in $\Cu(A^\sim)$. Since $[\pi(a)]=n$, we must
have $[\pi(a')]=0$ in $\Cu(\C)$, i.e., $[a']\in \Cu(A)$. Thus, $[a]-n[1]=[a']-0\cdot [1]$ in $\Cu^\sim(A)$. This shows that $\Cu(A)$
is mapped surjectively onto the positive elements of $\Cu^\sim(A)$.

Suppose that $\Cu(A^\sim)$ has cancellation of projections.
If $[a_1],[a_2]\in \Cu(A)$ are mapped into the same element of $\Cu^\sim(A)$ then
$[a_1]+k[1]=[a_2]+k[1]$ for some $k$, whence $[a_1]=[a_2]$.
\end{proof}

Among the axioms for the category $\CCu$ given in Subsection \ref{prelimsCu} we have not included that 0
be the smallest element of the ordered semigroups. This suits our purposes here since, unlike for $\Cu(A)$, $0$ may not the smallest element of $\Cu^\sim(A)$.
For example, $\Cu^\sim(\C)\cong \Z\cup\{\infty\}$, which is easily deduced from the alternative picture of $\Cu^\sim$ for unital C$^*$-algebras.

\begin{proposition}\label{CusimAinCu}
Let $A$ be a C$^*$-algebra of stable rank one. Then the ordered semigroup $\Cu^\sim(A)$ belongs to the category $\CCu$.
\end{proposition}

\begin{proof}
Let us first consider the case that $A$ is unital.
Notice, from the definition of $\Cu^\sim(A)$ for unital $A$, that the map $x\mapsto x+[1]$ 
(though clearly not additive) is an ordered set automorphism of $\Cu^\sim(A)$. Next, notice that for a finite set or increasing 
sequence $x_i\in \Cu^\sim(A)$, $i=1,2,\dots$, there is a sufficiently large $m\in\N$ such that $x_i+m[1]\geq 0$ for all $i$. 
But the ordered subsemigroup of positive elements of $\Cu^\sim(A)$ is isomorphic to $\Cu(A)$, by Lemma \ref{sur-inj}. So, the 
axioms of the category $\CCu$ may be routinely verified in $\Cu^\sim(A)$ by translating by a sufficiently large multiple of $[1]$ 
and using that $\Cu(A)$ is an ordered semigroup in  $\CCu$.

Let us now drop the assumption that $A$ is unital. Recall that the order on $\Cu^\sim(A)$ is the relative order as a subsemigroup of $\Cu^\sim(A^\sim)$. 
Suppose that $([a_i]-n_i[1])_{i=1}^\infty$
is an increasing sequence in $\Cu^\sim(A)$.  Then
\begin{align}\label{someseq}
0\leq [a_1]\leq [a_2]-(n_2-n_1)[1]\leq [a_2]-(n_3-n_1)[1]\leq \dots
\end{align}
is an increasing sequence in $\Cu(A^\sim)$ (here we identify $\Cu(A^\sim)$ with the positive elements of $\Cu^\sim(A^\sim)$ 
by Lemma \ref{sur-inj}). Let $[a]$ denote its  supremum. Since
$\Cu(\pi)$ applied to every term of \eqref{someseq} is equal to $n_1$ we have
$\Cu(\pi)([a])=n_1$. Thus $[a]-n_1[1]\in \Cu^\sim(A)$.
Notice that $[a]-n_1[1]$ is the supremum of $([a_i]-n_i[1])_{i=1}^\infty$ in $\Cu^\sim(A^\sim)$. Hence  $\Cu^\sim(A)$ is closed 
under the suprema of increasing sequences  inside $\Cu^\sim(A^\sim)$. This proves axiom O1 of the category $\CCu$
(see Subsection \ref{prelimsCu}).

If $[a]-n[1]\in \Cu^\sim(A)$ then $[(a-\epsilon)_+]-n[1]\in \Cu^\sim(A)$  for $0<\epsilon<\|a\|$
and
\[
[a]-n[1]=\sup_{\epsilon>0} ([(a-\epsilon)_+]-n[1]).\]
 This proves O2.
The   axioms O3 and O4  follow from their being valid  in $\Cu^\sim(A^\sim)$.
\end{proof}

Let $\phi\colon A\to B$ be a homomorphism and let $\phi^\sim\colon A^\sim\to B^\sim$ denote 
its unital extension. Then  $\Cu^\sim(\phi)\colon \Cu^\sim(A)\to \Cu^\sim(B)$ is defined as
\begin{align}
\Cu^\sim(\phi)([a]-n[1]):=[\phi^\sim(a)]-n[1].\label{Cusimphi}
\end{align}

\begin{remark}\label{Cusimphirem}
Let $\phi\colon A\to B$ be a homomorphism.

(i) If $A$ and $B$ are unital, and $\phi(1_A)=1_B$, then using the definition of $\Cu^\sim(A)$ and $\Cu^\sim(B)$ for unital C$^*$-algebras  the  map $\Cu^\sim(\phi)$ takes the form
\[
\Cu^\sim(\phi)([a]-n[1_A])=[\phi(a)]-n[1_B].
\]

(ii) If  $A$ is unital, without $B$ nor $\phi$ being necessarily unital, then using the definition of $\Cu^\sim(A)$  for unital C$^*$-algebras  the  map $\Cu^\sim(\phi)$ takes the form
\begin{align*}
\Cu^\sim(\phi)([a]-n[1_A]) &=[\phi(a)]-n[\phi(1_A)]\\
                                    &=[\phi(a)]+n[1_{B^\sim}-\phi(1_A)]-n[1_{B^\sim}].
\end{align*}

(iii) As with the functor $\Cu$, if $\phi$ is approximately unitarily equivalent to a homomorphism 
$\psi$  then $\Cu^\sim(\phi)=\Cu^\sim(\psi)$.
\end{remark}

\begin{proposition}\label{inCu}
Restricted to the C$^*$-algebras of stable rank one, $\Cu^\sim$ is a functor into the category $\CCu$ that preserves sequential inductive limits.
\end{proposition}
\begin{proof}
Let  $\phi\colon A\to B$ be a homomorphism, with $A$ and $B$ of stable rank one. Let us show that $\Cu^\sim(\phi)$ is a morphism in the category $\CCu$. Assume first that $\phi$ is unital. From Remark \ref{Cusimphirem} (i) we see that
$\Cu^\sim(\phi)$ is  additive. In particular, we have $\Cu^\sim(\phi)(x+[1])=\Cu^\sim(\phi)(x)+[1]$ for $x\in \Cu^\sim(A)$. In order to verify the properties of
a morphism of  $\CCu$ for the map $\Cu^\sim(\phi)$ we can use the same method used for the semigroup $\Cu^\sim(A)$: after translating by a sufficiently  large multiple of $[1]$, we reduce the verification of these properties to the restriction of $\Cu^\sim(\phi)$ to the subsemigroup of positive elements of $\Cu^\sim(A)$. By Lemma \ref{sur-inj}, this subsemigroup may  be identified with $\Cu(A)$. Finally, on $\Cu(A)$ the map $\Cu^\sim(\phi)$ coincides
with $\Cu(\phi)$, which is a morphism in  $\CCu$.

Let us drop the assumption that $A$ and $B$ are unital. 
From \eqref{Cusimphi} it is clear that $\Cu^\sim(A)$, viewed as a sub-semigroup of $\Cu^\sim(A^\sim)$, is preserved
by $\Cu^\sim(\phi^\sim)$. A routine verification (left to the reader) shows that the restriction of $\Cu^\sim(\phi^\sim)$ to
$\Cu^\sim(A )$ is a morphism of $\CCu$.

Finally, consider an inductive limit of stable rank one  C$^*$-algebras $A=\varinjlim A_i$. In order that $\Cu^\sim(A)$ be the inductive limit
of the $\Cu^\sim(A_i)$s  we must verify properties L1 and L2 of inductive limits in $\CCu$ (see Subsection \ref{prelimsCu}).
Unitizing in $A=\varinjlim A_i$ and taking
the Cuntz semigroup functor  we get
$\Cu(A^\sim)=\varinjlim \Cu(A_i^\sim)$.
Hence, for $[a]\in \Cu(A^\sim)$ there is an increasing sequence $([a_i])_{i=1}^\infty$ with supremum $[a]$, and where the $[a_i]$s come  from the ordered semigroups $\Cu(A_i^\sim)$. If $[\pi(a)]=n<\infty$, then $[\pi(a_i)]=n$ for $i$ large enough. Thus, $[a]-n[1]=\sup_i ([a_i]-n[1])$. This proves L1.

Let us prove L2. Let $[a_1]-n_1[1],[a_2]-n_2[1]\in \Cu^\sim(A_i)$ be such that the image of
$[a_1]-n_1[1]$ in $\Cu^\sim(A)$ is less than or equal to the image of $[a_2]-n_2[1]$. This is equivalent
to the image of  $[a_1]+n_2[1]$ in $\Cu(A^\sim)$ being less than or equal to the image of $[a_2]+n_2[1]$.
By the continuity of the functor $\Cu$, for every $\epsilon>0$ there is $j\geq i$ such that
the image of $[(a_1-\epsilon)_+]+n_2[1]$ in $\Cu(A_j^\sim)$ is less than or equal to the image of  $[a_2]+n_1[1]$
in $\Cu(A_j^\sim)$. That is, the image of $[(a_1-\epsilon)_+]-n_1[1]$
in $\Cu(A_j^\sim)$ is less than or equal to the image of  $[a_2]-n_2[1]$. This proves L2.
\end{proof}

\begin{proposition}\label{splitexact}
Let
\[
\xymatrix{
0\ar[r]& I \ar[r]^\phi\  & A\ar[r]^\psi & A/I\ar[r] &0
}
\]
be a short exact sequence of C$^*$-algebras, with $A$ of stable rank one.
We then have that

(i) $\mathrm{Im}(\Cu^\sim(\phi))=\mathrm{Ker}(\Cu^\sim(\psi))$;

(ii) if $\psi$ splits then $\Cu^\sim(\psi)$ is surjective and $\Cu^\sim(\phi)$ is injective.
\end{proposition}
\begin{proof}
(i) From $\psi\circ\phi=0$ we get that $\Cu^\sim(\psi)\circ\Cu^\sim(\phi)=0$. Therefore,
$\mathrm{Im}(\Cu^\sim(\phi))\subseteq \mathrm{Ker}(\Cu^\sim(\psi))$.

Let us prove that $\mathrm{Ker}(\Cu^\sim(\psi))\subseteq \mathrm{Im}(\Cu^\sim(\phi))$.
Let $[a]-n[1]\in \Cu^\sim(A)$ be  mapped to 0 by $\Cu^\sim(\psi)$. That is, $[\psi^\sim(a)]-n[1]=0$
in $\Cu((A/I)^\sim)$. By the cancellation of projections in
$\Cu((A/I)^\sim)$, this implies that $[\psi^\sim(a)]=[1_n]$ in $\Cu((A/I)^{\sim})$.
Since $\psi^\sim(a)$ is Cuntz equivalent to a projection, 0 is an isolated point of the spectrum of $\psi^\sim(a)$ 
(see \cite[Proposition 5.7]{brown-ciuperca}). So, for $f\in \mathrm C_0(\R^+)$ strictly positive and chosen suitably  $\psi^\sim(f(a))$  is a projection.
Rename $f(a)$ as $a$. Then $\psi^\sim(a)$ and  $1_n$ are Cuntz equivalent projections in $(A/I)^\sim\otimes\mathcal K$. 
Since $(A/I)^\sim\otimes\mathcal K$ is stably finite, $\psi^\sim(a)$ and $1_n$ are Murray-von
Neumann equivalent; furthermore, there exists  a unitary $u$ of the form $u=1+u'$, with $u'\in (A/I)^{\sim}\otimes\mathcal K$ such that $u^*\psi^\sim(a)u=1_n$. 
At this point, we go back to the start of the proof,  identify $A^\sim\otimes\mathcal K$
with $A^\sim\otimes \mathcal K\otimes \mathrm M_2$, and choose $a$ of the form 
\[
\begin{pmatrix}
a & 0\\
0 &0
\end{pmatrix}.
\] 
Notice then that the unitary $u$ may be chosen in 
$((A/I)^{\sim}\otimes\mathcal K\otimes e_{00})^\sim$. So we can lift 
\[
\begin{pmatrix}
u & 0\\
0 &u^*
\end{pmatrix}
\] 
to a unitary $v$ of the form $v=1+v'$, with $v'\in A^\sim\otimes\mathcal K\otimes \mathrm M_2$.  
Set $v^*av=a_1\in \mathrm M_2(A^\sim)\otimes \mathcal K$. Then $[a]=[a_1]$ and $\psi^\sim(a_1)=1_n$,
which says that $a_1\in I^\sim\otimes\mathcal K\otimes \mathrm M_2$. Thus, $[a]-n[1]=[a_1]-n[1]\in \mathrm{Im}(\Cu^\sim(\phi))$. 

(ii) Let $\lambda\colon A/I\to A$ be such that $\psi\circ \lambda=\mathrm{id}_{A/I}$.
Then $\Cu^\sim(\psi)\circ\Cu^\sim(\lambda)=\mathrm{id}_{\Cu^\sim(A/I)}$, which implies
that $\Cu^\sim(\psi)$ is surjective.

Let us prove that $\Cu^\sim(\phi)$ is injective.
Let $a,b\in I^\sim\otimes\mathcal K$ be positive elements such that $[\pi(a)]=[\pi(b)]=n$
and $[a]-n[1]=[b]-n[1]$ as elements of $\Cu^\sim(A)$. We want to show that $[a]-n[1]=[b]-n[1]$
in $\Cu^\sim(I)$. By the cancellation of projections in $\Cu(A^\sim)$ and $\Cu(I^\sim)$, we have that
$[a]=[b]$ in $\Cu(A^\sim)$ and we want to show that $[a]=[b]$ in $\Cu(I^\sim)$.

We may assume without loss of generality that $a=1_n+a'$ and  $b=1_n+b'$, for $a',b'\in I\otimes \mathcal K$. Let $\epsilon>0$. Since $a\cuntzle b$ in $A^\sim\otimes\mathcal K$, there exists $x\in A^\sim\otimes\mathcal K$ such that $(a-\epsilon/2)_+=x^*x$ and $xx^*\leq Mb$ for some $M>0$. Then by \cite[Theorem 5]{pedersen}, there exists 
a unitary $u\in (A^\sim\otimes\mathcal K)^\sim$ such that 
$u^*(a-\epsilon)_+u\leq Mb$. 

Recall that $\psi^\sim$ denotes the unital extension of $\psi$
to $A^\sim$, and that $\psi^\sim$ also denotes the homomorphism  $\psi^\sim\otimes \mathrm{id}$ extending $\psi^\sim$ to $A^\sim\otimes\mathcal K$. 
Let us furthermore continue to denote by $\psi^\sim$ -- rather than $(\psi^\sim)^\sim$ -- the unital extension of $\psi^\sim$  to $(A^\sim\otimes\mathcal K)^\sim$. We apply the same notational conventions to the homomorphism $\lambda$. Consider the elements
\begin{align*}
u_1 &=(\lambda^\sim\circ\psi^\sim)(u),\\
a_1 &=(uu_1^*)^*a(uu_1^*).
\end{align*}
We have $(\lambda^\sim\circ\psi^\sim)(uu_1^*)=1$, which implies that $uu_1^*$ is a unitary in $(I^\sim\otimes\mathcal K)^\sim$. Thus, $a\sim a_1$ in $I^\sim\otimes\mathcal K$. 
We also have that 
\[
(\lambda^\sim\circ\psi^\sim)(a_1)=(\lambda^\sim\circ\psi^\sim)(a)=1_n.
\] 
This implies that  $a_1=1_n+a_1'$ for $a_1'\in I\otimes\mathcal K$. 
From $u^*(a-\epsilon)_+u\leq Mb$ and the definition of $u_1$  
we get that $u_1^*(a_1-\epsilon)_+u_1\leq Mb$. 
Applying $\lambda^\sim\circ\psi^\sim$ on both
sides of this last inequality, and using that $(\lambda^\sim\circ\psi^\sim)(u_1)=u_1$, 
we get that $u_1^*1_nu_1=1_n$, i.e., $u_1$ commutes with $1_n$.

For $k=1,2,\dots$, let 
\[
v_k=1_n+(u_1-1_n)(e_k\otimes 1),
\]
where $e_k\in I^+$ is such that  
$(e_k\otimes 1)a_1'\to a_1'$ and $(e_k\otimes 1)u_1^*a_1'\to u_1^*a_1'$ when $k\to \infty$ (e.g., $(e_k)$ is an approximate unit for the hereditary subalgebra generated by the entries of $a_1'$ and $u_1^*a_1'$). Then $v_k$ belongs to $(I^\sim\otimes\mathcal K)^\sim$ for all $k$. We have
\[
v_ka_1'=1_na_1'+(u_1-1_n)(e_k\otimes 1)a_1'\to u_1a_1'.
\]
Here we have used that $(e_k\otimes 1)a_1'\to a_1'$.
Similarly, we deduce that $v_k^*a_1'\to u_1^*a_1'$. We also have that
\begin{align*}
v_k1_n &=1_n+(u_1-1_n)(e_k\otimes 1)1_n\\
& =1_n+1_n(u_1-1_n)(e_k\otimes 1)\\
&=1_nv_k,
\end{align*}
for all $k$ (here we have used that $u_1$ commutes with $1_n$). That is,  $v_k$ commutes with $1_n$ for all $k$. 
 
Since $I^\sim\otimes\mathcal K$ has stable rank one,  for each $k$ there exists a unitary $w_k$ in $(I^\sim\otimes\mathcal K)^\sim$ such that 
$\|v_k-w_k|v_k|\|<\frac{1}{k}$. Then $w_k^*a_1' w_k\to u_1^*a_1'u_1$ and 
$w_k^*1_n w_k\to 1_n$. It follows that 
$w_k^*a_1w_k\to u_1^*a_1u_1$. So $(a_1-\epsilon)_+\sim u_1^*(a_1-\epsilon)_+u_1$
in $I^\sim\otimes \mathcal K$. Hence,
\[
(a-\epsilon)_+\sim (a_1-\epsilon)_+\sim u_1^*(a_1-\epsilon)_+u_1\leq Mb,
\]
where the Cuntz comparisons are all taken in $I^\sim\otimes\mathcal K$.  Since $\epsilon>0$ is arbitrary, we conclude that $[a]\leq [b]$ in 
$\Cu(I^\sim)$. By symmetry, we also have $[b]\leq [a]$. So $[a]=[b]$ as elements of $\Cu(I^\sim)$.
\end{proof}

\begin{proposition}\label{isoembed} Let $A$ be a C$^*$-algebra of stable rank one.
The inclusion $A\hookrightarrow A\otimes\mathcal K$ (in the top corner) induces an isomorphism
of ordered semigroups
$\Cu^\sim(A)\to \Cu^\sim(A\otimes\mathcal K)$.
\end{proposition}

\begin{proof}
Consider the inductive limit
\[
A\hookrightarrow \mathrm M_2(A)\hookrightarrow \mathrm M_4(A)\hookrightarrow \dots \hookrightarrow A\otimes \mathcal K.
\]
By the continuity of $\Cu^\sim$ with respect to sequential inductive limits, it suffices to show that
the inclusion in the top corner $A\hookrightarrow \mathrm M_2(A)$ induces an isomorphism
at the level of $\Cu^\sim$.

Let us first assume that $A$ is unital. Then using the picture of $\Cu^\sim$ for unital C$^*$-algebras and  
homomorphisms with unital domain (see Remark \ref{Cusimphirem} (ii)) we see 
that the inclusion $A\hookrightarrow \mathrm M_2(A)$ induces an isomorphism in $\Cu^\sim$ if it induces  an isomorphism in $\Cu$. 
It is known that $A\hookrightarrow \mathrm M_2(A)$ induces an isomorphism in $\Cu$  (see \cite[Appendix 6]{coward-elliott-ivanescu}). So, the desired result follows for $A$ unital.

The non-unital case is reduced to the unital case as follows.
Let $A$ be a non-unital  C$^*$-algebra. Consider the diagram
\[
\xymatrix{
0\ar[r]& A \ar[r]\ar[d]  & A^\sim\ar[r]\ar[d]& \C \ar[r]\ar[d] &0\\
0\ar[r] & \mathrm M_2(A) \ar[r] & \mathrm M_2(A^\sim)\ar[r]& \mathrm M_2(\C) \ar[r] &0,
}
\]
where the rows form short exact sequences that split and the vertical
arrows are the natural inclusions.
Passing to the level of $\Cu^\sim$ we have
\[
\xymatrix{
0\ar[r]& \Cu^\sim(A) \ar[r]\ar[d]  & \Cu^\sim(A^\sim)\ar[r]\ar[d]& \Cu^\sim(\C) \ar[r]\ar[d] &0\\
0\ar[r] & \Cu^\sim(\mathrm M_2(A)) \ar[r] & \Cu^\sim(\mathrm M_2(A^\sim))\ar[r]& \Cu^\sim(\mathrm M_2(\C)) \ar[r] &0.
}
\]
A diagram chase -- as in the proof of the five lemma -- using the exactness of the rows of this diagram (in the sense of Proposition \ref{splitexact} (i) and (ii)),
and that the two rightmost vertical arrows  are isomorphisms, shows that $\Cu^\sim(A)\to \Cu^\sim(\mathrm M_2(A))$
is an isomorphism.
\end{proof}

\begin{proposition}\label{Cusimsums}
Let $A$ and $B$ be stable rank one C$^*$-algebras. Let $\iota_A\colon A\to A\oplus B$
and $\iota_B \colon B\to A\oplus B$ denote the standard inclusions. Then the map
$\gamma\colon \Cu^\sim(A)\oplus \Cu^\sim(B)\to\Cu^\sim(A\oplus B)$, given by
\[
\gamma(x+y):=\Cu^\sim(\iota_A)(x)+\Cu^\sim(\iota_B)(y)
\]
is an isomorphism of ordered semigroups.
\end{proposition}

\begin{proof}
Consider the diagram
\[
\xymatrix{
0\ar[r]& \Cu^\sim(A) \ar[r]^-{\alpha}\ar@{=}[d]  & \Cu^\sim(A)\oplus \Cu^\sim(B)\ar[r]^-{\beta}\ar[d]^\gamma & \Cu^\sim(B) \ar[r]\ar@{=}[d] &0\\
0\ar[r] & \Cu^\sim(A) \ar[r]^-{\Cu^\sim(\iota_A)} & \Cu^\sim(A\oplus B)\ar[r]^-{\Cu^\sim(\pi_B)}& \Cu^\sim(B) \ar[r] &0,
}
\]
where $\alpha(x)=(x,0)$, $\beta(x,y)=y$, and $\pi_B\colon A\oplus B\to B$ is the quotient map.
A simple diagram chase using the exactness of the rows of this diagram (in the sense of Proposition \ref{splitexact} (i) and (ii)) shows that $\gamma$
is an isomorphism, as desired.
\end{proof}

\subsection{Classification by $\Cu^\sim$: permanence properties}
Our main focus is the classification of homomorphisms from a given stable rank one C$^*$-algebra $A$ into
an arbitrary C$^*$-algebra of stable rank one. The classifying functor is   $\Cu^\sim$ together
with  the class $[s_A]\in \Cu^\sim(A)$ of a strictly positive element of $A$
(this class does not depend on the choice of the strictly positive element).  Since we will consider
repeatedly such a classification question, we introduce an abbreviated way of referring to it:

In the sequel, given a C$^*$-algebra $A$, we say  that \emph{the functor $\Cu^\sim$ classifies homomorphisms from $A$} if for any stable rank one C$^*$-algebra $B$, and any morphism
in $\CCu$
\[\alpha\colon \Cu^\sim(A)\to \Cu^\sim(B)\]
such that $\alpha([s_A])\leq [s_B]$, with $s_A$
and $s_B$ strictly positive elements of $A$ and $B$, there exists a homomorphism
$\phi\colon A\to B$ such that $\Cu^\sim(\phi)=\alpha$. Moreover, $\phi$ is unique up to approximate unitary equivalence (by unitaries in $B^\sim$).

\begin{lemma}\label{extendalpha}
Let $A$ and $B$ be C$^*$-algebras of stable rank one, with $B$ unital. Let $s_A\in A$
be a strictly positive element. Let $\alpha\colon \Cu^\sim(A)\to \Cu^\sim(B)$
be a  morphism in the
category $\CCu$ such that $\alpha([s_A])\leq [1]$. Then there exists a unique morphism in $\CCu$ $\tilde\alpha\colon \Cu^\sim(A^\sim)\to \Cu^\sim(B)$  that extends $\alpha$ and satisfies
$\tilde\alpha([1])=[1]$.
\end{lemma}
\begin{proof}
Let $\alpha$ be as in the statement of the lemma and suppose the extension $\tilde\alpha$ exists. For $[a]\in \W(A^\sim)$ such that
$[\pi(a)]=n<\infty$ we have, by linearity of $\tilde\alpha$, that
\begin{align}\label{alphatilde}
\tilde\alpha([a]-m[1])=\alpha([a]-n[1])+(n-m)[1].
\end{align}
Observe that the right side depends only on $\alpha$. Since the set
of elements $[a]-m[1]$, with $[a]\in \W(A^\sim)$ (i.e., $[\pi(a)]<\infty$), is dense in $\Cu^\sim(A^\sim)$,  the extension $\tilde\alpha$ of $\alpha$ is unique.

In order to prove the existence of $\tilde\alpha$, let \eqref{alphatilde}
stand for its definition   on the subsemigroup $\tilde W:=\{[a]-m[1]\mid [a]\in \W(A^\sim),m\in \Z\}$. It is clear from \eqref{alphatilde} that $\tilde\alpha$ is   additive
on this subsemigroup.
Let us prove that $\tilde\alpha$ preserves the order and far below relations on the elements of $\tilde W$, and that it preserves the supremum of increasing sequences in $\tilde W$ with supremum also in $\tilde W$.
Translating by a sufficiently large multiple of $[1]$, it suffices to prove these properties for the restriction of $\tilde\alpha$ to $\W(A^\sim)$ (see the proof of Proposition \ref{inCu}).

So let $[a]\in \W(A^\sim)$ and set $[\pi(a)]=n$. By  \cite[Theorem 1]{ciuperca-robert-santiago}, we have $n[1]\leq [a]+n[s_A]$, i.e.,
$0\leq [a]-n[1]+n[s_A]$ in $\Cu^\sim(A)$. Hence,
\begin{align}\label{alphapos}
0\leq \alpha([a]-n[1])+n\alpha([s_A])\leq \alpha([a]-n[1])+n[1]=\tilde\alpha([a]).
\end{align}
That is, $\tilde\alpha$ is positive.

Let $[a_1],[a_2]\in \W(A^\sim)$ be such that $[a_1]\leq [a_2]$. If $[\pi(a_1)]=[\pi(a_2)]$
then $\tilde\alpha([a_1])\leq \tilde\alpha([a_2])$ by \eqref{alphatilde} and the fact that $\alpha$
is order preserving. Suppose that  $[\pi(a_1)]<[\pi(a_2)]$. Let $\epsilon>0$.
By \cite[Lemma 7.1 (i)]{rordam-winter} there is $[c]\in \W(A^\sim)$ such that
\begin{align*}
[(a_1-\epsilon)_+]+[c] &\leq [a_2],\\
[\pi((a_1-\epsilon)_+)]+[\pi(c)] &=[\pi(a_2)].
\end{align*}
Thus,
\begin{align*}
\tilde\alpha([(a_1-\epsilon)_+])&\leq \tilde\alpha([(a_1-\epsilon)_+])+\tilde\alpha([c])\\
&= \tilde\alpha([(a_1-\epsilon)_+]+[c])\\
& \leq \tilde\alpha([a_2]).
\end{align*}
On the other hand, from \eqref{alphatilde} we deduce that
\[
\tilde\alpha([a])=\sup_{\epsilon>0}\tilde\alpha([(a-\epsilon)_+]),
\]
for $[a]\in \W(A^\sim)$. Thus, $\tilde\alpha([a_1])\leq \tilde\alpha([a_2])$, i.e., $\tilde\alpha$
is order preserving on $\W(A^\sim)$.

From \eqref{alphatilde} and the fact that $\alpha$ preserves the far below relation we get
$\tilde\alpha([(a-\epsilon)_+])\ll \tilde\alpha([a])$ for every $\epsilon>0$ and $[a]\in \W(A^\sim)$.
So if $[b]\ll [a]$ then $[b]\leq [(a-\epsilon)_+]\ll [a]$ for some $\epsilon>0$, and so
$\tilde\alpha([b])\leq \tilde([(a-\epsilon)_+])\ll \tilde \alpha([a])$. Thus, $\tilde\alpha$ preserves the
far below relation.

Finally, let $([a_i])$ be an increasing sequence with supremum $[a]\in \W(A^\sim)$.
Then $[\pi(a_i)]=[\pi(a)]<\infty$ for $i$ large enough. Now we deduce that $\tilde\alpha([a])=\sup_i \tilde\alpha([a_i])$  from \eqref{alphatilde} and the fact that $\alpha$
is supremum preserving.

We now extend $\tilde\alpha$ to $\Cu^\sim(A^\sim)$ by setting
\[
\tilde\alpha([a]-m[1]):=\sup_{\epsilon>0}\tilde\alpha([(a-\epsilon)_+]-m[1]).
\] 
(Notice that $[(a-\epsilon)_+]-m[1]\in \tilde W$ for all $\epsilon>0$.)  The subsemigroup $\tilde W$ is dense in $\Cu^\sim(A)$ (i.e., every element of $\Cu^\sim(A)$ is the supremum of a rapidly increasing sequence of elements of $\tilde W$) and 
belongs to the category $\mathbf{PreCu}$ (see \cite[Definition 2.1]{antoine-bosa-perera}). Thus, the properties of $\tilde\alpha$
already established on $\tilde W$ readily extend
to $\Cu^\sim(A^\sim)$ (see \cite[Theorem 3.3]{antoine-bosa-perera}).
\end{proof}

The meaning  of ``$\Cu$
classifies homomorphims from $A$" in part (i) of the following theorem is the same as the one defined above for $\Cu^\sim$, except with 
$\Cu$ in place of $\Cu^\sim$.

\begin{theorem}\label{permanence}
Let $A$  be a C$^*$-algebra of stable rank one.

(i) If $A$ is unital then the functor $\Cu^\sim$ classifies homomorphisms from $A$ if and only if $\Cu$
classifies homomorphisms from $A$.

(ii)  The functor $\Cu^\sim$ classifies homomorphisms from $A$ if and only if it classifies homomorphisms from $A^\sim$.

(iii) If $\Cu^\sim$ classifies homomorphisms from the  stable rank one C$^*$-algebras $A_i$, $i=1,2,\dots$, and $A=\varinjlim A_i$, then $\Cu^\sim$ classifies homomorphisms from $A$.

(iv) If $\Cu^\sim$ classifies homomorphisms from $A$ and $B$, with $B$ also of stable rank one, then it classifies homomorphisms from $A\oplus B$.

(v) If $\Cu^\sim$ classifies homomorphisms from $A$, then it classifies homomorphisms from $A'$ for any $A'$ stably isomorphic to $A$.
\end{theorem}

\begin{proof}
(i) Suppose that $\Cu^\sim$ classifies homomorphisms from $A$ and let us show that $\Cu$ classifies homomorphisms from $A$ too. 
Consider the uniqueness question first. Let $B$ be a stable rank one C$^*$-algebra and $\phi,\psi\colon A\to B$  homomorphisms
such that $\Cu(\phi)=\Cu(\psi)$. From $[\phi(1)]=[\psi(1)]$ and the fact that $\Cu(B^\sim)$ has cancellation of projections,
we get that $[1-\phi(1)]=[1-\psi(1)]$ in $\Cu(B^\sim)$. Now using the picture of $\Cu^\sim$ for homomorphism with unital domain
given in Remark \ref{Cusimphirem} (ii), we get that
\begin{align*}
\Cu^\sim(\phi)([a]-n[1]) &=[\phi(a)]+n[1-\phi(1)]-n[1]\\
&=[\psi(a)]+n[1-\psi(1)]-n[1]\\
&=\Cu^\sim(\psi)([a]-n[1]).
\end{align*}
That is, $\Cu^\sim(\phi)=\Cu^\sim(\psi)$. Since $\Cu^\sim$ classifies homomorphisms from $A$, we conclude that $\phi$ and $\psi$ are approximately unitarily equivalent,
as desired.

Let $\alpha\colon \Cu(A)\to \Cu(B)$ be a morphism in the category $\CCu$ such that $\alpha([1])\leq \alpha([s_B])$.
Since  $[\alpha(1)]$ is compact and $B$ has stable rank one, there exists a projection $p\in B\otimes\mathcal K$ such that 
$[\alpha(1)]=[p]$ (see \cite[Theorem 3.5]{brown-ciuperca}). Moreover, since $[p]\leq [s_B]$, $p$ may be chosen in $B$ 
(by Proposition \ref{sr1iso}).
Let us define $\tilde\alpha\colon \Cu^\sim(A)\to \Cu^\sim(B)$ by
\begin{align}\label{alphaextend}
\tilde\alpha([a]-n[1]):=\alpha([a])+n[1-p]-n[1].
\end{align}
It is easily verified that $\tilde\alpha$ is well defined, additive, and order preserving. 
Notice that the restriction of $\tilde\alpha$ to $\Cu(A)$ (identified with the subsemigroup of positive elements of $\Cu^\sim(A)$)
coincides with $\alpha$. Since $\alpha$ is a morphism in $\CCu$, it follows that $\tilde \alpha$ is a morphism in $\CCu$ too
(see the first paragraph of the proof of Proposition \ref{inCu}). Thus, there exists 
$\phi\colon A\to B$ such that $\tilde\alpha=\Cu^\sim(\phi)$. Restricting these morphisms to $\Cu(A)$ we get 
$\Cu(\phi)=\alpha$.

Suppose now that $\Cu$ classifies homomorphisms from $A$. Let $B$ be a C$^*$-algebra of stable rank one and let $\phi,\psi\colon A\to B$ be homomorphisms such that $\Cu^\sim(\phi)=\Cu^\sim(\psi)$. Then  identifying $\Cu(A)$ with the elements
greater than 0 of $\Cu^\sim(A)$ we arrive at $\Cu^\sim(\phi)|_{\Cu(A)}=\Cu^\sim(\psi)|_{\Cu(A)}$, i.e., $\Cu(\phi)=\Cu(\psi)$.
Since $\Cu$ classifies homomorphisms from $A$, the maps $\phi$ and $\psi$ are approximately unitarily equivalent.

Let $\alpha\colon \Cu^\sim(A)\to \Cu^\sim(B)$ be a morphism in the category $\CCu$ such that $\alpha([1])\leq [s_A]$.
Since $\Cu$ classifies homomorphisms from $A$, there exists $\phi\colon A\to B$ such that $\alpha|_{\Cu(A)}=\Cu(\phi)$.
Notice that $\alpha$ is uniquely determined by its restriction to $\Cu(A)$, by means of the formula 
\eqref{alphaextend}, and similarly for $\Cu^\sim(\phi)$. Thus, $\alpha=\Cu^\sim(\phi)$.

(ii) Suppose that $\Cu^\sim$ classifies homomorphisms from $A$. Let us prove that $\Cu$ classifies homomorphisms from $A^\sim$
(then by (i), $\Cu^\sim$ also classifies homomorphisms from $A^\sim$). 
We consider the uniqueness question first. 
Let $\phi,\psi\colon A^\sim\to B$ be homomorphisms such that $\Cu(\phi)=\Cu(\psi)$, with $B$ of stable rank one. 
Since $\phi(1)$ and $\psi(1)$ are Cuntz equivalent projections in a stable rank one C$^*$-algebra, they must be unitarily equivalent.
Thus, after conjugating by a unitary, we have $\phi(1)=\psi(1)=p$. From  $[\phi(a)]=[\psi(a)]$ for all $[a]\in \Cu(A^\sim)$
we deduce that $\Cu^\sim(\phi|_A)=\Cu^\sim(\psi|_A)$. By assumption $\Cu^\sim$ classifies homomorphisms from $A$.
Thus, $\phi|_A$ and $\psi|_A$ are approximately unitarily equivalent. By Proposition \ref{sr1hered} (i),
the unitaries implementing this equivalence may be chosen so that they commute with $p$. Hence, $\phi$ is approximately unitarily equivalent to $\psi$.

Let $\alpha\colon \Cu(A^\sim)\to \Cu(B)$ be a morphism in $\CCu$ such that $\alpha([1])\leq [s_B]$, with $B$ of stable rank one. 
Let us show that $\alpha$ is induced
by a homomorphism from $A^\sim$ to $B$. Since  $[\alpha(1)]\ll [\alpha(1)]$ and $B$ has stable rank one, there is a projection $p$ such that 
$[\alpha(1)]=[p]$ (see \cite[Theorem 3.5]{brown-ciuperca}). Moreover, since $[p]\leq [s_B]$, $p$ may be chosen in $B$ (by Proposition
\ref{sr1iso}). It is well known that the inclusion $pBp\hookrightarrow B$ induces an embedding $\Cu(pBp)\hookrightarrow \Cu(B)$.
The image of this embedding consists of the elements $[b]\in \Cu(B)$ such that $[p]$ is an order unit for $[(b-\epsilon)_+]$ for all
$\epsilon>0$. It follows that the image of $\alpha$ is contained in $\Cu(pBp)$ (viewed as a hereditary subsemigroup of $\Cu(B)$).
Thus, we may regard
$\alpha$ as a map from $\Cu(A^\sim)$ to $\Cu(pBp)$. The morphism $\alpha$ gives rise to a morphism
$\tilde \alpha\colon \Cu^\sim(A^\sim)\to \Cu^\sim(pBp)$ in the category $\CCu$ defined by $\tilde\alpha([a]-n[1]):=\alpha([a])-n[p]$.
By the classification of homomorphisms from $A$ by the functor
$\Cu^\sim$, there exists $\phi\colon A\to pBp$ such that $\Cu^\sim (\phi)=\tilde \alpha|_{\Cu^\sim(A)}$. 
Let $\phi^\sim\colon A^\sim\to pBp$ be the unital extension of $\phi$. 
We have $\Cu^\sim(\phi^\sim)=\tilde\alpha$ by  Lemma \ref{extendalpha}. Restricting these morphisms to $\Cu(A)$, we get $\Cu(\phi^\sim)=\alpha$.

Let us now assume that $\Cu^\sim$ classifies homomorphisms from $A^\sim$ and prove that it classifies homomorphisms from $A$.
Let us consider the uniqueness part first.
Let $\phi,\psi\colon A\to B$ be such that $\Cu^\sim(\phi)=\Cu^\sim(\psi)$, with $B$ of stable rank one.
By Lemma \ref{extendalpha}, we have $\Cu^\sim(\phi^\sim)=\Cu^\sim(\psi^\sim)$, where   $\phi^\sim$ and $\psi^\sim$ denote the unital extensions of $\phi$ and $\psi$ respectively.
 Since $\Cu^\sim$ classifies homomorphisms from $A^\sim$,  $\phi^\sim$
and $\psi^\sim$ are approximately unitarily equivalent. Hence, so are $\phi$ and $\psi$.

Let $\alpha\colon \Cu^\sim(A)\to \Cu^\sim(B)$ be a  morphism in $\CCu$
such that $\alpha([s_A])\leq [s_B]$. By Lemma \ref{extendalpha},
there is $\tilde\alpha\colon \Cu^\sim(A^\sim)\to \Cu^\sim(B^\sim)$ that extends $\alpha$ and
satisfies $\tilde\alpha([1])=[1]$. Let $\phi\colon A^\sim\to B^\sim$ be such that
$\Cu^\sim(\phi)=\tilde\alpha$.
Notice that $\phi$ must be unital, since $\Cu(\phi)([1])=[1]$. Hence $\phi(A)\subseteq B$ and
$\Cu^\sim(\phi|_A)=\alpha$.

(iii) This is a direct consequence of the continuity of the functor $\Cu^\sim$ with respect to sequential
inductive limits (see \cite{ciuperca-elliott-santiago} for a proof for the functor $\Cu$).

(iv) Let $\phi,\psi\colon A\oplus B\to C$ be homomorphisms that agree on $\Cu^\sim$, with $C$ of stable rank one.
Composing $\phi$ and $\psi$ with the inclusions $A\stackrel{\iota_A}{\hookrightarrow} A\oplus B$ and  $B\stackrel{\iota_B}{\hookrightarrow} A\oplus B$ and applying
the functor $\Cu^\sim$ we get  $\Cu^\sim(\phi|_A)=\Cu^\sim(\psi|_A)$ and $\Cu^\sim(\phi|_B)=\Cu^\sim(\psi|_B)$.
Since $\Cu^\sim$ classifies homomorphisms from $A$ and $B$, $\phi|_A$ and $\psi|_A$ are approximately unitarily equivalent, and $\phi|_B$ and $\psi|_B$ too.
Thus, $\phi$ and $\psi$ are approximately unitarily equivalent by Proposition \ref{sr1hered} (ii).

Let us prove the existence part of the classification. Let $\alpha\colon \Cu^\sim(A\oplus B)\to \Cu^\sim(C)$ be a morphism 
in the category $\CCu$ such that $\alpha([s_A+s_B])\leq [s_C]$. 
Since $\Cu^\sim$ classifies homomorphisms from $A$ and from $B$, there are homomorphisms
$\phi_A\colon A\to C$ and $\phi_B\colon B\to C$ that induce
$\alpha\circ \Cu^\sim(\iota_A)$ and $\alpha\circ \Cu^\sim(\iota_B)$ at the level of $\Cu^\sim$.
Moreover, using that $[\phi_A(s_A)]+[\phi_B(s_B)]\leq [s_C]$ we can choose $\phi_A$ and $\phi_B$ with orthogonal ranges. To achieve this, find $x\in \mathrm M_2(C)$  such that
$x^*x=\Big(\begin{smallmatrix}
\phi_A(s_A) &\\
& \phi_B(s_B)
\end{smallmatrix}\Big)$  and $xx^*\in C$.
(Here $C$ is identified with the top left corner of $\mathrm M_2(C)$).
Let $x=v|x|$ be the polar decomposition of $x$ in $\mathrm M_2(C)^{**}$. Then
\[
\phi_A'=v^*
\begin{pmatrix}
\phi_A &0\\
0&0
\end{pmatrix}v\hbox{ and }
\phi_B'=v^*
\begin{pmatrix}
 0&0\\
0&\phi_B
\end{pmatrix}v
\]
are the desired homomorphisms.

Let us now define $\phi\colon A\oplus B\to C$ by $\phi:=\phi_A\circ \pi_A+\phi_B\circ\pi_B$.
We have
\begin{align*}
\Cu^\sim(\phi)\circ \Cu^\sim(\iota_A) &=\Cu^\sim(\phi_A)=\alpha\circ \Cu^\sim(\iota_A),\\
\Cu^\sim(\phi)\circ \Cu^\sim(\iota_B) &=\Cu^\sim(\phi_B)=\alpha\circ \Cu^\sim(\iota_B).
\end{align*}
It follows from Proposition \ref{Cusimsums} that the ranges of $\Cu^\sim(\iota_A)$ and $\Cu^\sim(\iota_B)$ span $\Cu^\sim(A\oplus B)$.
Thus, $\Cu^\sim(\phi)=\alpha$.

(v) It suffices to show that the functor $\Cu^\sim$ classifies homomorphisms from $A$ if and only if it classifies homomorphisms from $A\otimes \mathcal K$.
Suppose that $\Cu^\sim$ classifies homomorphisms from $A$. Since $A\otimes\mathcal K=\varinjlim  \mathrm M_{2^i}(A)$, it is enough to show that $\Cu^\sim$
classifies homomorphisms from $\mathrm M_2(A)$ and apply part (ii).
The prove now proceeds as in \cite[Proposition 5 (ii)]{ers}. The properties of the functor $\Cu$
that are relevant to the arguments given in \cite[Proposition 5 (ii)]{ers} hold also for $\Cu^\sim$. Namely,
that the inclusion $A\hookrightarrow \mathrm M_2(A)$ induces an isomorphism at the level of $\Cu^\sim$, and that $\Cu^\sim$ induces the identity on approximately inner homomorphisms.

In order to prove that if  $\Cu^\sim$ classifies homomorphisms from $A\otimes\mathcal K$, then it
classifies homomorphisms from $A$ we proceed as in \cite[Proposition 5 (iv)]{ers}. Again,
the relevant properties of the functor $\Cu^\sim$ are that the inclusion $A\hookrightarrow A\otimes \mathcal K$ induces an isomorphism at the level of $\Cu^\sim$, and that $\Cu^\sim$ induces the identity on approximately inner homomorphisms.
\end{proof}

For C$^*$-algebras $A$ and $B$ of stable rank one, let us write $A\leftrightsquigarrow B$
if there is a sequence of stable rank one C$^*$-algebras $A=A_1,A_2,A_3,\dots,A_{n-1},A_n=B$ such that, for each $i$,
either $A_i$ is stably isomorphic to $A_{i+1}$, $A_i$ is the unitization of $A_{i+1}$, or
$A_{i+1}$ is the unitization of $A_i$. It is clear that $\leftrightsquigarrow$ is an equivalence relation. Notice that $A\leftrightsquigarrow B$ implies
$\mathrm K_1(A)\cong \mathrm K_1(B)$.
More importantly, we have the following corollary to the preceding theorem.

\begin{corollary}\label{squig}
If $A\leftrightsquigarrow B$ then
$\Cu^\sim$ classifies  homomorphisms from $A$ if and only if it classifies homomorphisms from $B$.
\end{corollary}

\subsection{Uniform continuity of the classification}
\begin{theorem}\label{uniform}
Let $A$ be a stable rank one C$^*$-algebra such that $\Cu^\sim$ classifies homomorphisms from $A$. Then for every
finite subset $F\subseteq A$ and $\epsilon>0$ there exists a finite subset $G\subseteq \Cu^\sim(A)$
such that for any two homomorphisms $\phi,\psi\colon A\to B$, with $B$ of stable rank one, if 
\begin{align}\label{uniformCulevel}
\begin{array}{c}
\Cu^\sim(\phi)(g')\leq \Cu^\sim(\psi)(g)\\
\Cu^\sim(\psi)(g')\leq \Cu^\sim(\phi)(g)
\end{array} 
\hbox{ for all $g',g\in G$ with $g'\ll g$,}
\end{align}
then there exists a unitary $u\in B^\sim$ such that
\begin{align}\label{uniformCstarlevel}
\|u^*\phi(f)u-\psi(f)\|<\epsilon\hbox{ for all }f\in F.
\end{align}
\end{theorem}

\begin{proof}
Suppose, by contradiction, that there is a pair $(F,\epsilon)$ such that for every
$G\subseteq \Cu^\sim(A)$ there exist  homomorphisms  $\phi_G\colon A\to B_G$
and $\psi_G\colon A\to B_G$
that satisfy \eqref{uniformCulevel}, but that do not satisfy \eqref{uniformCstarlevel}
for any unitary $u\in B_G^\sim$. Let us replace $B_G$ by $B_G^\sim$ and simply 
assume that $B_G$ is unital. (Notice that \eqref{uniformCulevel} continues to hold after doing this.) Let $B$ denote the quotient
of the C$^*$-algebra $\prod_G B_G$ of bounded nets $(b_G)$ by the ideal of nets such that $\|b_G\|\to 0$, where $G$ ranges through
the finite subsets of $\Cu^\sim(A)$. Let us show that $B$ has stable rank one. We will have this once we show that $\prod_G B_G$ has stable rank one.
Any given $b_G\in B_G$ can be approximated by elements of the form $u_G|b_G|$, with $u_G$ unitary. It follows that the elements of the form 
$(u_G)|(b_G)|$ (i.e., with polar decomposition) form a dense subset of  $\prod_G B_G$. But the elements with polar decomposition are in the closure of the invertible elements.
 Thus, $\prod_G B_G$ and $B$ have stable rank one.

Let $\phi,\psi\colon A\to B$
be the homomorphisms induced by $(\phi_G),(\psi_G)\colon A\to \prod_G B_G$ by passing to the quotient $B$.
We will show that on one hand $\phi$ and $\psi$ are not approximately unitarily equivalent, while on the other
$\Cu^\sim(\phi)=\Cu^\sim(\psi)$. Since $B$ is a stable rank one C$^*$-algebra, this contradicts
the assumption that $\Cu^\sim$ classifies homomorphisms from $A$.

The homomorphisms $\phi$ and $\psi$ cannot be approximately unitarily equivalent for if they 
were so, then there would be unitaries $(u_G)$ such that $\|u_G^*\phi_G(a)u_G-\psi_G(a)\|\to 0$.
But this would contradict our assumption that, for each $G$, $\phi_G$ and $\psi_G$ do not satisfy \eqref{uniformCstarlevel}.

Let $a\in \mathrm M_m(A^\sim)^+$ be such that $[\pi(a)]=n$.
For each $G$ that contains $[(a-\epsilon)_+]-n[1]$ and $[(a-\frac{\epsilon}{2})_+]-n[1]$ 
we have that 
\[
[(\phi_G^\sim(a)-\epsilon)_+)]-n[1]\leq [(\psi_G^\sim(a)-\frac{\epsilon}{2})_+]-n[1]. 
\]
Since $B_G$ is a C$^*$-algebra of stable rank one, $\Cu(B_G)$ has  cancellation
of projections. So, $[(\phi_G^\sim(a)-\epsilon)_+)]\leq [(\psi_G^\sim(a)-\frac{\epsilon}{2})_+]$. 
Thus, there exists $x_G\in \mathrm M_m(B_G)$ such that  
\[
(\phi_G^\sim(a)-2\epsilon)_+=x_G^*x_G\hbox{ and }h(\psi_G(a))\cdot x_Gx_G^*=x_Gx_G^*,
\]
where  $h\in \mathrm C_0(\R^+)^+$ is equal to 1 on $(\epsilon/2,\infty)$. Let us set $x$
equal to the image of  $(x_G)\in \prod_G \mathrm M_m(B_G)$ in $\mathrm M_m(B)$. Then $\phi^\sim((a-\epsilon)_+)=x^*x$ and $h(\psi^\sim(a))xx^*=xx^*$,
which implies that $[\phi^\sim((a-\epsilon)_+)]\leq [\psi^\sim(a)]$. Since $\epsilon>0$ is arbitrary, 
$[\phi^\sim(a)]\leq [\psi^\sim(a)]$ for all $a\in \mathrm M_m(A^\sim)^+$, and by the density of $\W(A^\sim)$ in $\Cu(A^\sim)$, we get
$\Cu^\sim(\phi)\leq \Cu^\sim(\psi)$. By symmetry, we also have  that $\Cu^\sim(\phi)\leq \Cu^\sim(\psi)$. Thus, $\Cu^\sim(\phi)=\Cu^\sim(\psi)$.  
\end{proof}

\section{Special cases of the classification}\label{specialcases}
In this section we discuss a few special cases of Theorem \ref{main}
in order to illustrate the argument used in the general case. (Thus, they will again 
be dealt with when we prove Theorem \ref{main} in full generality.)

Our point of departure is Ciuperca and Elliott's result \cite[Theorem 4.1]{ciuperca-elliott} 
that (in the present sense) the functor $\Cu$ classifies homomorphisms from $\mathrm C_0(0,1]$ and $\mathrm C[0,1]$. 
By Theorem \ref{permanence} (i) and (ii),  we then have that $\Cu^\sim$
also classifies homomorphisms from $\mathrm C_0(0,1]$ and $\mathrm C[0,1]$.

\subsection{Trees}
By a  compact tree let us understand a  finite, connected,  1-dimensional simplicial complex without loops.

\begin{theorem}\label{trees}
Let $T$ be either a compact tree, or a compact tree with a point removed, or a finite disjoint union of such spaces. Then the functor $\Cu^\sim$ classifies homomorphisms from $\mathrm C_0(T)$.
\end{theorem}

\begin{proof}
By Theorem \ref{permanence} (iv), it suffices to consider only trees (with or without a point removed).
The proof proceeds by induction on the number of edges, the base case being the trees
$(0,1]$ and $[0,1]$. Suppose that the theorem is true for all compact trees with less than $n$ edges, and for all such spaces with  one  point removed. Let $T$ be a compact tree with $n$ edges.
Let $x$ be an endpoint of $T$ and consider the space $T_0$
obtained deleting the edge that contains $x$ from $T$. Then $T_0$ is a finite disjoint union of trees with less than $n$ edges and an endpoint removed. 
By Theorem \ref{permanence} (iv) and the induction hypothesis, $\Cu^\sim$ classifies homomorphisms from $\mathrm C_0(T_0)$.
Since $\mathrm C(T)\cong (\mathrm C_0(T_0)\oplus \mathrm C_0(0,1])^\sim$, the functor $\Cu^\sim$ classifies  homomorphisms
from $\mathrm C(T)$ by Theorem \ref{permanence} (ii) and (iv). Finally, since  $\mathrm C(T)$ is the unitization of $\mathrm C_0(T\backslash\{y\})$ for $y\in T$, $\Cu^\sim$ also classifies homomorphisms from $\mathrm C_0(T\backslash\{y\})$ for any $y\in T$. This completes the induction.
\end{proof}

\subsection{The algebra $q\C$.}
Let us show that the functor $\Cu^\sim$ classifies homomorphisms from the C$^*$-algebra
\[
q\C=\left \{\, f\in \mathrm M_2(\mathrm C_0(0,1]) \mid f(1)=
\begin{pmatrix}
* & 0\\
0 & *
\end{pmatrix}
\right \}.
\]
This algebra is the unitization of
\[
\left \{\, f\in \mathrm M_2(\mathrm C_0(0,1]) \mid f(1)=
\begin{pmatrix}
* & 0\\
0 & 0
\end{pmatrix}
\right \},
\]
which, in turn, is stably  isomorphic to $\mathrm C_0(0,1]$. Thus, $q\C \leftrightsquigarrow \mathrm C_0(0,1]$, and so $\Cu^\sim$ classifies homomorphisms from $q\C$ by Corollary \ref{squig}.

\subsection{Razak's building blocks.}
In \cite{razak}, Razak proves a classification result for simple inductive limits
of building blocks of the form $\mathrm M_m(\C)\otimes R_{1,n}$, where
\[
R_{1,n}=\left \{f\in \mathrm \mathrm M_n(\mathrm C[0,1])\mid
\begin{array}{l}
f(0)=\lambda 1_{n-1}\\
 f(1)=\lambda 1_n
\end{array},
\hbox{ for some }\lambda\in \C
\, \right \}.
\]
Let us see that $\Cu^\sim$ classifies homomorphisms from $R_{1,n}$ (whence, also
from the inductive limits of algebras of the form $\mathrm M_m(\C)\otimes R_{1,n}$). We have that
$R_{1,n}^\sim$ (i.e., the unitization of $R_{1,n}$) is the subalgebra of $\mathrm M_n(\mathrm C[0,1])$ of functions such that
\[
f(0)=
\begin{pmatrix}
\lambda 1_{n-1} &         \\
                         & \mu
\end{pmatrix} \hbox{ and }
f(1)=\lambda 1_n,
\]
for some $\lambda\in \C$.
Now notice that $R^\sim_{1,n}$ is also the unitization of
\[
\left\{\,
f\in \mathrm M_n(\mathrm C_0[0,1)) \mid f(0)=
\begin{pmatrix}
0_{n-1} & 0\\
0 & \mu
\end{pmatrix}
\, \right\}.
\]
This last algebra is stably isomorphic to $\mathrm C_0[0,1)$. Thus $R_{1,n} \leftrightsquigarrow \mathrm C_0(0,1]$, and so the functor $\Cu^\sim$ classifies homomorphisms from $R_{1,n}$.

\subsection{Prime dimension drop algebras.}
Let $p$ and $q$ be relatively prime, with $q>p$. Consider the 1-dimensional NCCW complex
\[
Z_{p,q}=\{\, f\in \mathrm M_{pq}(\mathrm C[0,1])\mid f(0)\in 1_q\otimes \mathrm M_p,\, f(1)\in 1_p\otimes \mathrm M_q \,\}.
\]
By Remark \ref{stablyNCCW} (iii), $Z_{p,q}$ is stably isomorphic to
\[
A_{p,q}=
\left \{\,
f\in \mathrm M_q(\mathrm C[0,1])\mid
\begin{array}{l}
f(0)=\mu 1_{p}\\
 f(1)=\lambda 1_q
\end{array},
\hbox{ for some }\lambda,\mu\in \C
\,\right \}.
\]
Let us show that $\Cu^\sim$ classifies homomorphisms from $A_{p,q}$, and hence also from $Z_{p,q}$.
We have $A_{p,q}\leftrightsquigarrow A_{q-p,p}$, since $A_{p,q}$ and $A_{q-p,q}$ both have the same unitization.
Let us assume  that $2p<q$ (otherwise, passing to $A_{q-p,q}$ we have $2(q-p)<q$ ). Then $A_{p,q}$ is isomorphic to a full hereditary subalgebra of $A_{p,q}'$, where $A_{p,q}'$ is composed of the functions $f\in \mathrm M_q(\mathrm C[0,1])$ such that
\[
f(0)=
\begin{pmatrix}
\Lambda\cdot 1_{p} &\\
                      & 0_{q-2p}
\end{pmatrix}
\hbox{ and }
f(1)=\lambda 1_q,
\]
with $\lambda\in \C$ and $\Lambda\in \mathrm M_2(\C)$. Thus $A_{p,q}\leftrightsquigarrow A_{p,q}'$. The unitization of $A_{p,q}'$ changes only the fiber at 0, and is also the unitization
of
\[
A_{p,q}''=
\left \{\,
f\in \mathrm M_q(\mathrm C_0[0,1))\mid
f(0)=
\begin{pmatrix}
\Lambda\cdot 1_{p} &\\
                      & \mu 1_{q-2p}
\end{pmatrix},\, \Lambda\in \mathrm M_2(\C),\mu\in \C
\,\right \}.
\]
Thus $A_{p,q}'\leftrightsquigarrow A_{p,q}''$. By Remark \ref{stablyNCCW} (iii), the algebra $A_{p,q}''$, is in turn
stably isomorphic to
\[
A_{p,q}'''=
\left \{\,
f\in \mathrm M_{q-p}(\mathrm C_0[0,1))\mid
f(0)=
\begin{pmatrix}
\lambda\cdot 1_{p} &\\
                      & \mu 1_{q-2p}
\end{pmatrix},\, \lambda,\mu\in \C
\,\right \}.
\]
Notice finally that  the unitization of $A_{p,q}'''$ is also the unitization of  $A_{p,q-p}$. Therefore, $A_{p,q}\leftrightsquigarrow A_{p,q-p}$. Since $p$ and $q$ are relatively prime, the continuation of this process leads to an algebra of the form $A_{1,d}$. The algebra $A_{1,d}$  is a a full hereditary subalgebra of $\mathrm M_d(\mathrm C_0[0,1))^\sim$.
It follows that   $A_{p,q} \leftrightsquigarrow \mathrm C_0(0,1]$, and so $\Cu^\sim$ classifies homomorphisms from $A_{p,q}$ for $p$ and $q$ relatively prime.

\section{General case of the classification}\label{proofofmain}
Throughout this section we will assume the notation for NCCW complexes introduced  in
Subsection \ref{NCCWprelims}.

\subsection{Reduction lemmas}
Recall that the 1-dimensional NCCW complex defined
in \eqref{nccw} is  
determined (up to isomorphism) by the data $((e_j)_{j=1}^l,(f_i)_{i=1}^k,Z^{\phi_0},Z^{\phi_1})$.
Here $(e_j)_{j=1}^l$ and $(f_i)_{i=1}^k$ are vectors of natural numbers,  $Z^{\phi_0}$ and $Z^{\phi_1}$
are $k\times l$ matrices of non-negative integers, and they are all related by the condition stated in Remark \ref{stablyNCCW} (ii).
As a matter of convenience, we will use the following notations:  $\mathbf e :=(e_j)_{j=1}^l$, $\mathbf f:=(f_i)_{i=1}^k$, and  for each $t=0,1$, the  columns of  $Z^{\phi_t}$ will be  denoted by   $\mathbf c^{\phi_t}_1,\mathbf c^{\phi_t}_2,\dots,\mathbf c^{\phi_t}_l$. 

Let $A$ be the 1-dimensional NCCW complex determined by the data $(\mathbf e,\mathbf f,Z^{\phi_0},Z^{\phi_1})$. 
Let us examine the effect that has on these data
to pass from $A$ to $A'$, where $A'$ is a 1-dimensional NCCW complex that is either stably isomorphic to $A$, isomorphic to the unitization of $A$, or
the result of ``removing a unit" from $A$ (i.e., $A$ is isomorphic to the unitization of $A'$). 
The verification of the following claims is straightforward (in each case,  $(\tilde {\mathbf e},\tilde{\mathbf f}, Z^{\tilde \phi_0},Z^{\tilde \phi_1})$
denotes the data associated to $A'$):
\begin{enumerate}
\item
Let $Z^{\tilde \phi_0}=Z^{\phi_0}$, $Z^{\tilde \phi_1}=Z^{\phi_1}$,  and choose the vectors $\tilde {\mathbf e}$, $\tilde{\mathbf f}$  arbitrarily (while still satisfying
the condition of Remark \ref{stablyNCCW} (ii)). Then $A'$ is stably isomorphic to $A$ (by Remark \ref{stablyNCCW} (iii)). 

\item
Let $Z^{\tilde \phi_0}=(Z^{\phi_0},\mathbf c_{l+1}^{\tilde \phi_0})$ and $Z^{\tilde \phi_1}=(Z^{\phi_1},\mathbf c_{l+1}^{\tilde \phi_1})$, where 
  $\mathbf c_{l+1}^{\tilde \phi_0}$ and $\mathbf c_{l+1}^{\tilde \phi_1}$ are defined by 
\[
Z^{\phi_0}\mathbf e+\mathbf c_{l+1}^{\tilde \phi_0}=Z^{\phi_1}\mathbf e+\mathbf c_{l+1}^{\tilde \phi_1}=\mathbf f.
\] 
Let $\tilde{\mathbf e}=(\mathbf e,1)$ and $\tilde {\mathbf f}=\mathbf f$. (That is, new columns are appended to the matrices
$Z^{\phi_0}$ and $Z^{\phi_1}$ so that the resulting homomorphisms $\tilde \phi_0$ and $\tilde \phi_1$ be unital, and a new entry is appended to the vector
$\mathbf e$ with the value of 1.) Then  $A'$ is isomorphic to the unitization of $A$. 

\item
Suppose that $Z^{\phi_0}\mathbf e=Z^{\phi_1}\mathbf e=\mathbf f$
and that  $e_{j_0}=1$ for some index $1\le j_0\leq l$. 
Then the previous transformation can be  reversed   as follows:
Let $Z^{\tilde \phi_0}$ and $Z^{\tilde \phi_1}$ be the matrices obtained by deleting
the $j_0$-th column from $Z^{\phi_0}$ and $Z^{\phi_1}$ respectively. Let $\tilde{\mathbf e}$ be the vector obtained by removing the $j_0$-th entry from $\mathbf e$, 
and let  $\tilde {\mathbf f}=\mathbf f$.  Then $A'$  is the result of ``removing a unit" from $A$. That is, the unitization of $A'$ is isomorphic to $A$.
\end{enumerate}

Observe that permuting the columns or rows of the matrices $Z^{\phi_0}$ and $Z^{\phi_1}$
does not change the isomorphism class of the resulting 1-dimensional NCCW complex, as long as the 
same permutation is simultaneously done on both matrices. We can also interchange the $i$-th rows of $Z^{\phi_0}$ and $Z^{\phi_1}$ while leaving
the other rows unchanged.

Let us say that the 1-dimensional NCCW complex given by the data $(\mathbf e,\mathbf f,Z^{\phi_0},Z^{\phi_1})$, is in \emph{reduced form} if
\begin{enumerate}
\item [(R1)] $\mathbf e=(1,1,\dots,1)$,

\item [(R2)] $Z^{\phi_0}\mathbf e=Z^{\phi_1}\mathbf e=\mathbf f$ (i.e., the maps $\phi_0$ and $\phi_1$ are unital),

\item [(R3)] $Z^{\phi_0}_{i,j}=0$ or $Z^{\phi_1}_{i,j}=0$ for all $i=1,2,\dots,l$ and $j=1,2,\dots,k$.
\end{enumerate}

\begin{lemma}\label{reduced}
Every 1-dimensional NCCW complex is equivalent, by the relation $\leftrightsquigarrow$, to one in reduced form.
\end{lemma}
\begin{proof}
By setting $\tilde {\mathbf e}=(1,1,\dots,1)$ without changing $Z^{\phi_0}$, $Z^{\phi_1}$ and $\mathbf f$, we get a  C$^*$-algebra stably isomorphic to the one that we started with. If it does not satisfy (R2), we add a unit to it. 
In this way, we get a C$^*$-algebra equivalent to the original one and satisfying (R1) and (R2).

To get the property (R3), let us  proceed as follows: Let us assume that $A$ already satisfies (R1) and (R2). Let $j$ be an index between $1$ and $l$. We perform the following steps (the input  data of every step is denoted by $(\mathbf e,\mathbf f,Z^{\phi_0},Z^{\phi_1})$ and the output data by
$(\tilde {\mathbf e},\tilde{\mathbf f}, Z^{\tilde \phi_0},Z^{\tilde \phi_1})$):

\emph{Step 1}. Remove the unit corresponding to $e_j=1$.
By the remarks made above, this step amounts to deleting the $j$-th columns of $Z^{\phi_0}$ and  $Z^{\phi_1}$ and the $j$-th entry of $\mathbf e$. 

\emph{Step 2}. Define   $\tilde {\mathbf f}$ by
\begin{align*}
\tilde f_i=\max\Big(\sum_{j=1}^l Z_{i,j}^{\phi_0},\sum_{j=1}^l Z_{i,j}^{\phi_1}\Big),
\end{align*}
for all $i=1,2,\dots,k$, while leaving the rest of the data unchanged. By the remarks made above, the resulting algebra is stably isomorphic to the  algebra obtained in the previous step. 

\emph{Step 3}.  Add a unit to the algebra obtained in the previous step. By the remarks made above, this amounts to inserting new columns (at the $j$-th spot) to the matrices
$Z^{\phi_0}$ and $Z^{\phi_1}$. The matrices $Z^{\tilde \phi_0}$ and $Z^{\tilde \phi_1}$ of the new
algebra agree with the corresponding matrices of the original algebra, except for their $j$-th columns, which now satisfy that $Z^{\tilde \phi_0}_{i,j}=0$ or $Z^{\tilde\phi_1}_{i,j}=0$ for all $i=1,2,\dots,k$. 

Repeating  Steps 1-3 for each $j=1,2,\dots, l$ we get an algebra that satisfies (R1), (R2), and (R3) and is  equivalent to the original algebra.
\end{proof}

A 1-dimensional NCCW complex in reduced form can be completely recovered from
the matrix $Z^\phi=Z^{\phi_0}-Z^{\phi_1}$. Indeed, the entries of $Z^{\phi_0}$ and $Z^{\phi_1}$ are recovered from $Z^\phi$ using (R3). The vector $\mathbf f$ is then recovered using  (R2) (while $\mathbf e$ is fixed by (R1)). We shall refer to $Z^\phi$ as the matrix associated to the 1-dimensional NCCW complex. Let us denote the columns of $Z^{\phi}$
by $\mathbf c^\phi_j$, with $j=1,2,\dots,l$.  Notice that any matrix $Z^\phi$ with integer entries and with  rows adding up to 0 is associated to a 1-dimensional NCCW complex in reduced form. 

\begin{lemma}\label{addsubtract}
Let $A$ be a 1-dimensional NCCW complex in reduced form and with associated matrix $Z^\phi$. Let $1\leq j_1,j_2,j_3\leq l$ be distinct indices. Then $A\leftrightsquigarrow A'$, where $A'$ is the 1-dimensional NCCW complex in reduced form with associated  matrix $Z^{\tilde \phi}$ equal to $Z^\phi$, except for the $j_2$-th and $j_3$-th columns, which are given by  $\mathbf c_{j_2}^{\tilde \phi}=\mathbf c_{j_2}^\phi-\mathbf c_{j_1}^\phi$ and $\mathbf c_{j_3}^{\tilde \phi}=\mathbf c_{j_3}^\phi+\mathbf c_{j_1}^\phi$.
\end{lemma}

\begin{proof}
Let us describe a sequence of steps going from $A$ to $A'$ (as in the proof of the previous lemma, the input  data of every step is denoted by $(\mathbf e,\mathbf f,Z^{\phi_0},Z^{\phi_1})$ and the output data by
$(\tilde {\mathbf e},\tilde{\mathbf f}, Z^{\tilde \phi_0},Z^{\tilde \phi_1})$): 

\emph{Step 1}.
Remove the unit corresponding to $e_{j_2}=1$. The output matrix $Z^{\tilde \phi}$ is obtained by deleting the $j_2$-th column of $Z^{\phi}$.

\emph{Step 2}.
Pass to the stably isomorphic algebra with the same matrices $Z^{\phi_0}$ and  $Z^{\phi_1}$, $\tilde {\mathbf e}=(1,1,\dots,2,1,\dots,1)$, where $\tilde e_{j_1}=2$, and
vector $\tilde {\mathbf f}$ given by
\begin{align*} 
\tilde f_i=\max(2Z_{i,j_1}^{\phi_0}+\sum_{j\neq j_1}  Z_{i,j}^{\phi_0},
2Z_{i,j_1}^{\phi_1}+\sum_{j\neq j_1} Z_{i,j}^{\phi_1}),
\end{align*}
for all $i=1,2,\dots,k$.

\emph{Step 3}. Add a unit, inserting the new column at the $j_2$-th spot. The resulting associated matrix $Z^{\tilde \phi}$ agrees with the initial one except at  the 
$j_2$-th column. Since $Z^{\tilde \phi}\tilde {\mathbf e}=0$, and $\tilde e_{j_1}=2$, we get that the $j_2$-th column of $Z^{\tilde\phi}$ must be equal to $\mathbf c_{j_2}^\phi-\mathbf c_{j_1}^\phi$. Observe that the current algebra is not in reduced form, since
$\tilde e_{j_1}=2$.

\emph{Step 4}.
Remove the unit corresponding to $e_{j_3}=1$. This has the effect of deleting the $j_3$-th column from the matrix of the previous step.

\emph{Step 5}. Pass to the stably isomorphic algebra with  $\tilde {\mathbf{e}}=(1,1,\dots,1)$, the same matrices $Z^{\phi_0}$ and $Z^{\phi_1}$, and
vector $\tilde{\textbf f}$  given by
\begin{align*}
\tilde f_i=\max(\sum_{j=1}^l Z_{i,j}^{\phi_0},\sum_{j=1}^l Z_{i,j}^{\phi_1}),
\end{align*}
for all $i=1,2,\dots,k$.

\emph{Step 6}. Add a unit, inserting the new column at the $j_3$-th spot. A straightforward computation yields that the new column is equal to $\mathbf c_{j_3}^\phi-\mathbf c_{j_1}^\phi$.
\end{proof}

\subsection{Proof of  Theorem \ref{main}}

The following lemma is a straightforward consequence of the Mayer-Vietoris
sequence in K-theory applied to the pull-back diagram \eqref{nccw} (see also 
\cite{elliott2}).

\begin{lemma}\label{K10}
Let $A$ be a 1-dimensional NCCW complex as in \eqref{nccw}.
Then $\mathrm K_1(A)=0$ if and only if  $\mathrm K_0(\phi_0)-\mathrm K_0(\phi_1)$ is surjective.
\end{lemma}

\begin{proposition}\label{reduction}
Let $A$ be a 1-dimensional NCCW complex such that $\mathrm K_1(A)=0$. Then $A\leftrightsquigarrow\mathrm C[0,1]$.
\end{proposition}
\begin{proof}
By Lemma \ref{reduced}, we may assume that $A$ is in reduced form. Let  $Z^\phi=Z^{\phi_0}-Z^{\phi_1}$ be the associated matrix of $A$.
By the previous lemma, $Z^{\phi}$ defines a surjection from $\Z^l$ to $\Z^k$. 
Let us first find a suitable reduction of the first row of $Z^\phi$. Suppose  that the first row of $Z^\phi$ has at least three non-zero entries. Using Lemma \ref{addsubtract}, we can reduce the smallest absolute value among  these three entries. We may continue doing this until there are at most two non-zero entries in the first row of $Z^\phi$. We may assume (after performing some permutations) that these entries are  $Z_{1,1}^\phi$ and $Z_{1,2}^\phi$. Since the rows of $Z^\phi$ add up to zero, we  
must have $Z^{\phi}_{1,1}=-Z^{\phi}_{1,2}$. Furthermore, since  $Z^\phi$ is surjective, we must have $Z^{\phi}_{1,2}=\pm 1$.
Without loss of generality, let us assume that $Z^{\phi}_{1,2}=-1$. (To multiply a row of $Z^{\phi}$ by $-1$, interchange the corresponding rows of $Z^{\phi_0}$ and $Z^{\phi_1}$; this does not change the isomorphism class of the 1-dimensional NCCW complex.) The first row of the resulting matrix $Z^\phi$ is equal to   $(1,-1,0,\dots,0)$.

Now consider the second row of $Z^{\phi}$. Observe that at least one entry $Z_{2,j}^\phi$ for $j\geq 3$ must be non-zero, by the surjectivity of $Z^\phi$ (since a 2$\times$2 matrix with rows that add up to 0 cannot be surjective). If there exist two entries $Z_{2,j}^\phi,Z_{2,j'}^\phi$ with $j,j'\geq 3$ that are non-zero, we can reduce the smallest absolute value among these two entries by suitable applying  Lemma \ref{addsubtract} (for the indices $2,j,j'$). We may continue doing this until there is exactly one non-zero entry $Z_{2,j}^{\phi}$, with $j\geq 3$.
Let us assume without loss of generality that this entry is $Z_{2,3}^{\phi}$. Since the rows of
$Z^{\phi}$ add up to zero, the  matrix obtained by deleting the first column of $Z^\phi$ defines a surjection from $\Z^k$ to $\Z^l$. Furthermore, deleting all  but the first two rows of this matrix results in the 2$\times $2  matrix
\[
\begin{pmatrix}
-1 & 0\\
Z_{2,2}^\phi & Z_{2,3}^\phi
\end{pmatrix},
\] 
which is also surjective, from $\Z^2$ to $\Z^2$.
Thus, we must have $Z_{2,3}^\phi=\pm 1$. Let us assume without loss of generality
that $Z_{2,3}^\phi=-1$. We can apply Lemma \ref{addsubtract} to the first three columns of $Z^\phi$  (with $j_1=3$), to get that $Z_{2,1}^\phi=0$. Then, we  must have that $Z_{2,2}^\phi=1$, since the rows of $Z^\phi$ add up to zero.
As a result, we get
\[
Z^\phi=
\begin{pmatrix}
1 & -1 & 0 & 0 & \cdots\\
0 & 1 & -1 & 0 & \cdots\\
\cdots & \cdots & \cdots & \cdots &\cdots \\
\cdots & \cdots & \cdots & \cdots &\cdots 
\end{pmatrix}.
\]
Continuing this process with the third, fourth  row of $Z^{\phi}$, etc,  until we have considered all rows, we get 
\[
Z^\phi=
\begin{pmatrix}
1 & -1 & 0 & 0 & 0 & \cdots\\
0 & 1 & -1 & 0 & 0 &\cdots\\
0 & 0 & 1 & -1 &  0 &\cdots \\
\cdots & \cdots & \cdots & \cdots &\cdots\\
0 &\cdots & 0 & 1 & -1 & \cdots  
\end{pmatrix}.
\]
The 1-dimensional NCCW complex in reduced form with this associated matrix is isomorphic to $\mathrm C[0,1]\oplus \C^{k-l-1}$. Removing units successively for each of the $\C$ summands, we get $\mathrm C[0,1]$.
\end{proof}

\begin{proof}[Proof of Theorem \ref{main}]
By Theorem \ref{permanence}, it suffices to show that $\Cu^{\sim}$ classifies homomorphisms from a 1-dimensional NCCW complex with trivial $\mathrm K_1$-group. Furthermore, by Proposition \ref{reduction} and Corollary \ref{squig}, it suffices to show that $\Cu^\sim$ classifies homomorphisms from $\mathrm C[0,1]$. But as pointed out at the beginning Section \ref{specialcases}, this is essentially Ciuperca and Elliott's \cite[Theorem 1]{ciuperca-elliott} (combined with Theorem \ref{permanence} (i)).
\end{proof}

\begin{corollary}\label{isomorphism}
Let $A$ and $B$ be inductive limits of 1-dimensional NCCW complexes with trivial $\mathrm K_1$-group. Then
$A\otimes\mathcal K\cong B\otimes\mathcal K$ if and only if $\Cu^\sim(A)\cong \Cu^\sim(B)$. If the isomorphism 
from $\Cu^\sim(A)$ to $\Cu^\sim(B)$ maps $[s_A]$ into $[s_B]$, where $s_A$ and $s_B$ are strictly positive elements
of $A$ and $B$ respectively, then $A\cong B$. Moreover,  in this case the isomorphism from $\Cu^\sim(A)$
to $\Cu^\sim(B)$ lifts to an isomorphism from $A$ to $B$.
\end{corollary}

\begin{proof}
See the proof of \cite[Corollary 1]{ciuperca-elliott-santiago}.
\end{proof}

\section{Simple C$^*$-algebras}\label{simple}
Here we show that if a C$^*$-algebra $A$ is simple and an inductive limit of 1-dimensional NCCW complexes then $\Cu^\sim(A)$ is determined 
by $\mathrm K_0(A)$  and the cone of traces of $A$
(and their pairing). In the case that $A\otimes\mathcal K$ contains a projection, this is
essentially a consequence of Winter's $\mathcal Z$-stability result \cite{winter} for simple C$^*$-algebras of finite decomposition
rank, and of the computation of $\Cu(A)$ obtained by Brown and Toms in \cite{brown-toms} for simple unital $\mathcal Z$-stable $A$ 
($\mathcal Z$ denotes the Jiang-Su algebra).
If $A\otimes\mathcal K$ is projectionless, a different route in the computation of $\Cu^\sim(A)$ must be followed, since  there is currently no version of Winter's result
for  projectionless C$^*$-algebras (although such a result is likely to be true). So, in this case we rely 
on more  ad hoc methods  to compute $\Cu(A)$, which we show to be determined solely by the cone of traces of $A$. 
We then compute $\Cu^\sim(A)$ in terms of
$\mathrm K_0(A)$, the cone of traces of $A$, and their pairing.

Let $A$ be a  C$^*$-algebra. Let $\T_{0}(A)$ denote the cone of densely finite, positive, lower semicontinuous traces on $A$. The cone $\T_0(A)$ is endowed with the topology of pointwise
convergence on elements of the Pedersen ideal of $A$. 
We shall consider various  spaces of functions on $\T_0(A)$:
\begin{align*}
\mathrm{Aff}_+(\T_0(A)) & :=\left \{f\colon \T_0(A)\to [0,\infty)\mid
\begin{array}{l}
f\hbox{ is linear, continuous, and}\\
f(\tau)>0\hbox{ for }\tau\neq 0
\end{array}
\right\},
\\
\mathrm{LAff}_+(\T_0(A)) & :=\{f\colon \T_0(A)\to [0,\infty]\mid
\exists (f_n) \hbox{ with }f_n\uparrow f\hbox{ and }f_n\in \mathrm{Aff}_+(\T_0(A))\},\\
\mathrm{LAff}_+^\sim(\T_0(A))  &:=\left \{f\colon \T_0(A)\to (-\infty,\infty]\mid f=f_1-f_2,
\begin{array}{l}
 f_1\in \mathrm{LAff}_+(\T_0(A)),\\
f_2\in \mathrm{Aff}_+(\T_0(A))\end{array}
\right\}.
\end{align*}

Since lower semicontinuous traces on $A$ extend (uniquely) to 
$A\otimes\mathcal K$, positive elements  $a\in (A\otimes\mathcal K)^+$, and Cuntz classes $[a]\in \Cu(A)$, give rise to functions on $\T_0(A)$:
\begin{align*}
\widehat {a}(\tau) &:=\tau(a),\hbox{ for }\tau\in \T_0(A),\\
\widehat{[a]} &: =\sup_n \widehat{(a^{\frac 1 n})}.
\end{align*}

\begin{remark}
The following  facts are either known or easily verified.

(i) If $a$ is a positive element in the Pedersen ideal of $A\otimes\mathcal K$   then $\widehat a\in \mathrm{Aff}_+(\T_0(A))$. 

(ii) If $a$ is an arbitrary positive element  in $A\otimes\mathcal K$ then  $\widehat a\in \mathrm{LAff}_+(\T_0(A))$. (Because $\widehat a=\sup_n \widehat{ (a-\frac 1 n)_+ }$ and
$(a-\frac 1 n)_+$ belongs to the Pedersen ideal of $A\otimes\mathcal K$ for all $n>0$.) 

(iii) If $A$ is simple then
the range of the map $a\mapsto \widehat a$, with $a\in (A\otimes\mathcal K)^+$, is exactly the space $\mathrm{LAff}_+(\T_0(A))$ (see \cite[Remarks 5.14 and 6.9]{ers}).
\end{remark}

\subsection{Case with projections}
Let $A$ be a simple inductive limit of 1-dimensional NCCW complexes. Assume, furthermore, that $A\otimes\mathcal K$ contains at least one
non-zero projection $r$. Then $r(A\otimes\mathcal K)r$ is a full hereditary subalgebra of $A\otimes \mathcal K$ -- by the simplicity of $A\otimes \mathcal K$. Thus, by Brown's theorem,  $A$ and $r(A\otimes\mathcal K)r$ are stably isomorphic and therefore
have isomorphic $\Cu^\sim$-ordered semigroups. Since $r(A\otimes\mathcal K)r$ is unital, we may use
the picture of $\Cu^\sim$ for unital C$^*$-algebras. Thus, a general element of 
$\Cu^\sim(A)$ has the form $[a]-[q]$, where 
$[a]\in \Cu(r(A\otimes\mathcal K)r)=\Cu(A)$ and $q\in A\otimes\mathcal K$ is a projection.

Consider the set $\mathrm K_0(A)\sqcup \mathrm{LAff}_{+}^\sim(\T_0(A))$. 
Let us define on this set an ordered semigroup structure. On the subsets $\mathrm K_0(A)$ and $\mathrm{LAff}_{+}^\sim(\T_0(A))$ the addition operation agrees with the addition
with which these sets are endowed. For mixed sums, let us define
\[
([p]-[q])+\alpha=\widehat{[p]}-\widehat{[q]} +\alpha\in \mathrm{LAff}_{+}^\sim(\T_0(A)),
\]
where $[p]-[q]\in \mathrm{K}_0(A)$ and $\alpha\in \mathrm{LAff}_{+}^\sim(\T_0(A))$.
The order again restricts to the natural order on the subsets  $\mathrm K_0(A)$ and $\mathrm{LAff}_{+}^\sim(\T_0(A))$.
For $[p]-[q]$ and $\alpha$ as before, let $\alpha\leq [p]-[q]$ if
 $\alpha\leq \widehat{[p]}-\widehat{[q]}$,
and  $[p]-[q]\leq \alpha$ if $\widehat{[p]}-\widehat{[q]} < \alpha$.

Let us define a map from $\Cu^\sim(A)$ to $\mathrm K_0(A)\sqcup \mathrm{LAff}_{+}^\sim(\T_0(A))$ by
\begin{align}\label{isoCuLTAp}
[a]-[q] \mapsto 
\left\{
\begin{array}{ll}
[p]-[q] & \hbox{if $[a]=[p]$, with $p$ a projection in $A\otimes\mathcal K$}\\
\widehat{[a]}-\widehat{[q]} & \hbox{otherwise.}
\end{array}
\right.
\end{align}

\begin{proposition}\label{Cusimprojections}
Let $A$ be a simple inductive limit of 1-dimensional NCCW complexes. Suppose that $A\otimes\mathcal K$ contains at least one non-zero projection. Then the 
map \eqref{isoCuLTAp} is an isomorphism of ordered semigroups.
\end{proposition}

\begin{proof}
Let $r$ be a non-zero projection in $A\otimes\mathcal K$. 
Then $r(A\otimes\mathcal K)r$ is simple and  unital. Since the decomposition rank is well behaved with 
respect to inductive limits and hereditary subalgebras, $r(A\otimes\mathcal K)r$ has   decomposition rank 
at most 1 (see \cite[Section 3.3 and Proposition 3.10]{kirchberg-winter}).  It follows by Winter's \cite[Theorem 5.1]{winter} that  
$r(A\otimes\mathcal K)r$ absorbs tensorially the Jiang-Su algebra $\mathcal Z$. But $A$ and $r(A\otimes\mathcal K)r$ are stably
isomorphic. Thus, $A$  is $\mathcal Z$-absorbing too (by \cite[Corollary 3.1]{toms-winter}). It now follows from the computation of 
$\Cu(A)$ in \cite[Theorem 2.5]{brown-toms}  that  
\begin{align*}
[a] \mapsto 
\left\{
\begin{array}{ll}
[p]& \hbox{if $[a]=[p]$, with $p$ a projection in $A\otimes\mathcal K$}\\
\widehat{[a]}& \hbox{otherwise,}
\end{array}
\right.
\end{align*}
is an ordered semigroup isomorphism from $\Cu(A)$ to $\mathrm V(A)\sqcup \mathrm{LAff}_{+}(\T_0(A))$.
Here $\mathrm V(A)$ denotes the semigroup of Murray-von Neumann classes of  projections
of $A\otimes \mathcal K$. Recall that we view  $\Cu^\sim(A)$ as the semigroup of formal differences $[a]-[q]$,
 with $[a]\in \Cu(A)$ and $q\in A\otimes\mathcal K$ a projection. A straightforward calculation then shows that \eqref{isoCuLTAp} 
is also an isomorphism of ordered semigroups.
\end{proof}

\subsection{Projectionless case}
Let us now turn to the computation of $\Cu^\sim(A)$ in the stably projectionless case.

\begin{proposition}\label{Cuprojectionless}
Let $A$ be a simple inductive limit of 1-dimensional NCCW complexes. Suppose that
$A\otimes\mathcal K$ is projectionless.
Then the mapping $[a]\mapsto \widehat{[a]}$, from $\Cu(A)$ to
$\mathrm{LAff}_{+}(\T_0(A))$, is an isomorphism of ordered semigroups.
\end{proposition}

\begin{proof}
Let $A_i$ be 1-dimensional NCCW complexes such that $A=\varinjlim (A_i,\phi_{i,j})$.
By \cite[Theorem 4.6]{toms-dyn}, the Cuntz semigroup of a 1-dimensional NCCW complex has strict comparison. That is,
if $\widehat{[a]}\leq (1-\epsilon)\widehat{[b]}$ for some $\epsilon>0$, then $[a]\leq [b]$.
Since strict comparison passes  to inductive limits, $\Cu(A)$ has strict comparison too. It is known that for a simple C$^*$-algebra with strict comparison, the map $[a]\mapsto \widehat{[a]}$ is injective on the
complement of the subsemigroup of Cuntz classes $[p]$, with $p$ a projection (see \cite{ers}, \cite{brown-toms}). Since in our case $A\otimes \mathcal K$ is projectionless, we conclude that $[a]\mapsto \widehat{[a]}$ is injective.

In order to prove surjectivity, it suffices to show that for each $[a]\in \Cu(A)$ and $\lambda\in \Q^+$
there exists $[b]$ such that $\lambda\widehat{[a]}=\widehat{[b]}$ (see \cite[Corollary 5.8]{ers}). In fact, it suffices to show
this for $\lambda=\frac{1}{n}$, $n\in \N$. Using that $\Cu(A)$ has strict comparison, this can be reduced
to proving the following property of almost divisibility: 
\begin{enumerate}
\item[(D)] For all $x\in \Cu(A)$, $x'\ll x$, and $n\in \N$,
there exists $y\in \Cu(A)$ such that $n\widehat{y}\leq \widehat{x}$ and 
$\widehat{x'}\leq (n+1)\widehat{y}$. 
\end{enumerate}
The argument to prove that $\Cu(A)$ has this property runs along similar lines as the proof of \cite[Theorem 3.4]{toms-rigid}. 
We sketch it here briefly: First, notice that in order to prove (D)  it suffices
to verify it for $x$ belonging to a dense subset of $\Cu(A)$, as it then extends easily to all $x$. (Recall that by dense subset we mean one that
for every $x$ there is a rapidly increasing sequence of elements in the given subset with supremum $x$.)
For $B\subseteq (A\otimes\mathcal K)^+$ dense and closed under functional calculus, the set $\{[b]\mid b\in B\}$
is dense in $\Cu(A)$ (this is a consequence of R\o rdam's \cite[Proposition 4.4]{rordam}).  
Thus, we may assume that $x=[\phi_{i,\infty}(a)]$, with $a\in (A_i\otimes\mathcal K)^+$ for some $i$.  Let $0<\epsilon<\|a\|$.
Since $A$ is simple (and non-type I)  there exists $j$ such that $\phi_{i,j}(a)$
and $\phi_{i,j}((a-\epsilon)_+)$ have rank at least $n+1$ on every irreducible representation of $A_j\otimes\mathcal K$. Rename
$\phi_{i,j}(a)$ as $a$.

\emph{Claim.} There exists a positive element $b\in A_j\otimes\mathcal K$ such that 
\begin{align}\label{rankab}
n\cdot \rank_\pi b &\leq \rank_\pi a,\\
\rank_\pi (a-\epsilon)_+ &\leq (n+1)\cdot \rank_\pi b,\label{rankab2}
\end{align}
for every irreducible representation $\pi$ of $A_j\otimes\mathcal K$. 

\emph{Proof of claim.} In order to find $b$, we use 
that $A_j$ is a 1-dimensional NCCW complex. Let $A_j$ be given by the pull-back diagram 
\eqref{nccw}. By \cite[Theorem 3.1]{antoine-perera-santiago}, we may identify $\Cu(A_j)$ with the ordered semigroup of pairs $((n_i)_{i=1}^k,f)$, with $n_i\in \{0,1,\dots,\infty\}$ and
$f\in \mathrm{Lsc}([0,1],\{0,1,\dots,\infty\}^l)$ such that 
\[
Z^{\phi_0}(n_i)\leq f(0)\hbox{ and }Z^{\phi_1}(n_i)\leq f(1).
\]
Say $[a]$ corresponds to the pair $((n_i)_{i=1}^k,f)$. Then we can choose $[b]$ as the pair $((\lfloor\frac{n_i}{n}\rfloor)_{i=1}^k,\lfloor\frac{f}{n}\rfloor)$. A simple calculation shows that
$[b]$ satisfies \eqref{rankab} and \eqref{rankab2}, as desired.

Let $[b]$ be as in the previous claim. Then \eqref{rankab} and \eqref{rankab2} imply that $n\widehat{[b]}\leq \widehat{[a]}$ 
and $\widehat{[(a-\epsilon)_+]}\leq (n+1)\widehat{[b]}$. Moving $[a]$ and $[b]$ forward to $\Cu(A)$ we get (D).
\end{proof}

Let us now define a pairing between elements of $\Cu^\sim(A)$ and traces in $\T_0(A)$. That is, we define a map from $\Cu^\sim(A)$ to $\mathrm{LAff}_{+}^\sim(\T_0(A))$. 
Let $A$  be a simple, projectionless inductive limit of 1-dimensional NCCW complexes. 
In the sequel, we  identify $\Cu(A)$ with $\mathrm{LAff}_{+}(\T_0(A))$ by the isomorphism
given in Proposition \ref{Cuprojectionless}. 
Let $[a]\in \Cu(A^\sim)$ be such that $[\pi(a)]=n<\infty$. Since $[a]$ and $n[1]$ are mapped to the same element in the quotient by $A$,
there exists $[b]\in \Cu(A)$ such that $n[1]\leq [a]+[b]$ (see \cite[Proposition 1]{ciuperca-robert-santiago}). Since we are identifying
$\Cu(A)$ with $\mathrm{LAff}_{+}(\T_0(A))$, we write $n[1]\leq [a]+\beta$
with $\beta \in \mathrm{LAff}_{+}(\T_0(A))$. Recall that every function of $\mathrm{LAff}_{+}(\T_0(A))$ is the supremum of 
an increasing sequence of functions in $\mathrm{Aff}_{+}(\T_0(A))$. So, we may choose $\beta$  in $\mathrm{Aff}_{+}(\T_0(A))$ such that $n[1]\leq [a]+\beta$. By Lemma \ref{algebraicorder}, $n[1]$ sits as a summand
of any element above it. Thus,  there exists $\gamma\in \mathrm{LAff}_{+}(\T_0(A))$ such that
\begin{align}\label{betagamma}
n[1]+\gamma=[a]+\beta.
\end{align}
In summary, for every $[a]\in \Cu(A^\sim)$ such that $[\pi(a)]=n<\infty$ there exist $\beta\in \mathrm{Aff}_{+}(\T_0(A))$
and $\gamma\in \mathrm{LAff}_{+}(\T_0(A))$  such that \eqref{betagamma} holds.
\begin{lemma}Assume the notation and hypotheses of the preceding paragraph.
The assignment
\begin{align}\label{mapintoLTA}
[a]-n[1]\mapsto \gamma-\beta.
\end{align}
is a well defined map from $\Cu^\sim(A)$ to $\mathrm{LAff}_{+}^\sim(\T_0(A))$.
\end{lemma}

\begin{proof}
Let $[a_1]-n_1[1]$ and $[a_2]-n_2[1]$ be elements of $\Cu^\sim(A)$ such that $[a_1]-n_1[1]\leq [a_2]-n_2[1]$.
Let $\beta_1\in \mathrm{Aff}_{+}(\T_0(A))$ and 
$\gamma_1\in \mathrm{LAff}_{+}(\T_0(A))$ be functions such that \eqref{betagamma} holds for $[a_1]-n_1[1]$ and define $\beta_2$
and $\gamma_2$ similarly for $[a_2]-n_2[1]$. 
The following computations are performed in
$\Cu^\sim(A)$:
\begin{align*}
[a_1]-n_1[1] &\leq [a_2]-n_2[1]\\
[a_1]-n_1[1]+\beta_1+\beta_2 &\leq [a_2]-n_2[1]+\beta_1+\beta_2\\
n_1[1]+\gamma_1-n_1[1]+\beta_2 &\leq n_2[1]+\gamma_2-n_2[1]+\beta_1\\
\gamma_1+\beta_2 &\leq \gamma_2+\beta_1.
\end{align*}
That is,  $\gamma_1-\beta_1\leq \gamma_2-\beta_2$.
This shows at once that the map \eqref{mapintoLTA}
is well defined and order preserving.
\end{proof}

For $[a]-n[1]\in \Cu^\sim(A)$ let us denote by $([a]-n[1])^{\widehat{}}$
the function $\gamma-\beta$, with $\beta$ and $\gamma$ as in \eqref{betagamma}.
We have shown in the previous lemma that this map is well defined and order preserving.
It is also clear from its definition that it is additive.
We can now put an order and an additive structure on $\mathrm K_0(A)\sqcup \mathrm{LAff}_{+}^\sim(\T_0(A))$ just as we did before 
in the case that $A\otimes\mathcal K$
contained a projection. This time we use the map $[a]-n[1]\mapsto ([a]-n[1])^{\widehat{}}$ to define
mixed sums and the order relation. Moreover, we can define a map from $\Cu^\sim(A)$ to $\mathrm K_0(A)\sqcup \mathrm{LAff}_{+}^\sim(\T_0(A))$ by
\begin{align}\label{isoCuLTA}
[a]-n[1] \mapsto
\left\{
\begin{array}{ll}
[p]-n[1] & \hbox{if $[a]=[p]$, with $p$ a projection in $A^\sim\otimes\mathcal K$}\\
([a]-n[1])^{\widehat{}} & \hbox{otherwise.}
\end{array}
\right.
\end{align}

\begin{proposition} \label{Cusimprojectionless}
Let $A$ be a simple C$^*$-algebra such that $A\otimes \mathcal K$ is projectionless. Suppose that $A$ is either an inductive limit of 1-dimensional NCCW complexes or that $A\otimes \mathcal Z\cong A$, where $\mathcal Z$ denotes the Jiang-Su algebra. Then the map \eqref{isoCuLTA}  is an isomorphism of ordered semigroups.
\end{proposition}

\begin{proof}
Let us assume that $A$ is an inductive limit of 1-dimensional NCCW complexes.  
We will prove that \eqref{isoCuLTA} defines a bijection from  $\Cu^\sim(A)$ to the ordered
semigroup $\mathrm K_0(A)\sqcup \mathrm{LAff}_{+}^\sim(\T_0(A))$. The verification
that this map is additive and order preserving is left to the reader. (Notice that we have already established
that its restrictions to $\mathrm{K}_0(A)$ and $\mathrm{LAff}_{+}^\sim(\T_0(A))$ are ordered semigroup
maps.)

Since $A^\sim$ is stably finite, the Cuntz equivalence of projections in $A^\sim\otimes\mathcal K$
amounts  to their Murray-von Neumann equivalence. Thus, the restriction of the map  \eqref{isoCuLTA} to the subsemigroup of $\Cu^\sim(A)$ of elements of the form $[p]-n[1]$, with $p$ a projection,  is a bijection with $\mathrm K_0(A)$.

Let us show that the restriction of the map \eqref{isoCuLTA} to the elements $[a]-n[1]$, with
$[a]\neq [p]$ for any projection $p$, is a bijection with   $\mathrm{LAff}_{+}^\sim(\T_0(A))$.

Let us show injectivity. Let $[a_1]-n_1[1]$ and $[a_2]-n_2[1]$ be elements of $\Cu^\sim(A)$ such that 
$([a_1]-n_1[1])^{\widehat{}}=([a_2]-n_2[1])^{\widehat{}}$.
Choose $\beta_1\in \mathrm{Aff}_{+}(\T_0(A))$ and 
$\gamma_1\in \mathrm{LAff}_{+}(\T_0(A))$ such that $([a_1]-n_1[1])^{\widehat{}}=\gamma_1-\beta_1$ and find $\beta_2$ and $\gamma_2$ 
similarly for $[a_2]-n_2[1]$. Then  
\begin{equation}\label{abetagamma}
[a_1]+\beta=n[1]+\gamma=[a_2]+\beta,
\end{equation}
with $\beta=\beta_1+\beta_2$ and $\gamma=\gamma_2+\beta_1$.

\emph{Claim}. If $[a_1]\neq [p]$, for any projection $p$, then $[(a_1-\epsilon)_+] + \beta\ll [a_1]+\beta$
for any $0<\epsilon<\|a_1\|$.

\emph{Proof of claim}: Let $g_\epsilon\in \mathrm C_0(0,\|a_1\|]^+$ be a function with support $(0,\epsilon]$. We have
\[
[(a_1-\frac{\epsilon}{2})_+] +[g_\epsilon(a_1)]+ \beta\leq [a_1]+\beta.
\]
We cannot have $g_\epsilon(a_1)= 0$; otherwise  0 would be an isolated point of the spectrum of $a_1$ 
and this would imply $[a_1]=[p]$ for some projection $p$. Since $[g_\epsilon(a_1)]\in \Cu(A)$, we view it
as a function in $\mathrm{LAff}_{+}(\T_0(A))$. Then, there is $\delta>0$ such that 
$(1+\delta)\beta\leq [g_\epsilon(a_1)]+ \beta$. Since $\beta\ll (1+\delta)\beta$ (because $\beta$ is continuous), we get that $\beta\ll [g_\epsilon(a_1)]+ \beta$. So,
\[
[(a_1-\epsilon)_+] + \beta\ll [(a_1-\frac{\epsilon}{2})_+] +[g_\epsilon(a_1)]+ \beta\leq [a_1]+\beta.
\]
This proves the claim.

Now from \eqref{abetagamma} we deduce that $[(a_1-\epsilon)_+]+\beta\ll [a_1]+\beta=[a_2]+\beta$
for any $\epsilon>0$. By weak cancellation, this implies $[(a_1-\epsilon)_+]\ll [a_2]$ for any $\epsilon>0$.
Hence $[a_1]\leq [a_2]$, and by symmetry, $[a_2]\leq [a_1]$. This shows that \eqref{mapintoLTA}
is injective on the elements $[a]-n[1]$ with $[a]\neq [p]$, for any projection $p\in A^\sim \otimes \mathcal K$.

Let us show surjectivity onto $\mathrm{LAff}_{+}^\sim(\T_0(A))$. Since the functions in $\mathrm{LAff}_{+}(\T_0(A))$
are attainable by elements in $\Cu(A)$,  it is enough to
show that $-\gamma$, with $\gamma\in  \mathrm{Aff}_{+}(\T_0(A))$, is attainable. 
For every $\epsilon>0$ we have $\gamma\ll (1+\epsilon)\gamma$. Thus, there exists  
a positive contraction $a$ in $\mathrm M_m(A)$ for some $m$, and a number $\delta>0$, such that 
\[
\gamma\leq \widehat{[(a-\delta)_+]}\leq \widehat{[a]}\leq  (1+\epsilon)\gamma.
\]
Let $c_\delta\in \mathrm C_0[0,1)$ be positive on $[0,\delta)$ and $0$ elsewhere. Let
$c_\delta(a)\in \mathrm M_m(A^\sim)$ denote the positive element obtained applying functional calculus
on $a$ and $1_m$. Then
\[
[c_\delta(a)]+[(a-\delta)_+]\leq m[1]\leq [c_\delta(a)]+[a].
\]
Passing to $\Cu^\sim(A)$ we have
\[
 [c_\delta(a)]-m[1]+[(a-\delta)_+]\leq 0\leq [c_\delta(a)]-m[1]+[a].
\]
From this we deduce that $([c_\delta(a)]-m[1])^{\widehat{}}\leq -\widehat{[a-\delta]_+}\leq -\gamma$.
Similarly, we have $-(1+\epsilon)\gamma\leq ([c_\delta(a)]-m[1])^{\widehat{}}$.
In summary, we have found $c_\delta(a)\in \mathrm M_m(A^\sim)$ such that
\[
-(1+\epsilon)\gamma\leq ([c_\delta(a)]-m[1])^{\widehat{}}\leq -\gamma.
\]
Notice that the spectrum of $a$ has no gaps, since $A\otimes\mathcal K$ is projectionless.
Thus, $[c_\delta(a)]$ is not the class of a projection in $A^\sim\otimes\mathcal K$.
In order to attain the function $-\gamma$ exactly, we consider the increasing sequence of functions $-(1+\frac{1}{n})\gamma$,
$n=1,2,\dots$, with pointwise supremum $-\gamma$. By the previous discussion, between any two consecutive terms of this sequence 
we can interpolate an element of the form $([c_n]-m_n[1])^{\widehat{}}$, where $[c_n]$ is not the class
of a projection. While proving the injectivity of the map \eqref{isoCuLTA} before, we have in fact
shown that if $([a_1]-n_1[1])^{\widehat{}}\leq ([a_2]-n_2[1])^{\widehat{}}$, and $[a_1]$ is not the class
of a projection, then $[a_1]-n_1[1]\leq [a_2]-n_2[1]$. Thus, we have that the sequence $[c_n]-m_n[1]$, $n=1,2,\dots$, is itself increasing.
It follows that $([c]-m[1])^{\widehat{}}=-\gamma$, where $[c]-m[1]$ the supremum of the sequence $[c_n]-m_n[1]$, $n=1,2,\dots$.

The proof above applies equally well to the case that $A\otimes \mathcal Z\cong A$, since all that was needed was that
the conclusion of Proposition \ref{Cuprojectionless} be valid (together with the fact that $A$ is simple and stably projectionless). In the case that $A\otimes\mathcal Z\cong A$, the conclusion of Proposition \ref{Cuprojectionless}
follows from \cite[Theorem 6.6]{ers}. 
\end{proof}

In the following corollary, we combine the computations of $\Cu^\sim(A)$ for both cases,
with and without projections, with the classification result of the previous section.
The pairing between $\mathrm K_0(A)$ and $\T_0(A)$ alluded to in the statement of this corollary is defined 
as follows: In the case that $A\otimes\mathcal K$ contains a projection,
the pairing between $\mathrm K_0(A)$ and $\T_0(A)$ is given by
\[
([p]-[q],\tau)\mapsto \tau(p)-\tau(q),
\]
while, in the stably projectionless case, the pairing is obtained by restricting the pairing 
between $\Cu^\sim(A)$ and $\T_0(A)$ to $\mathrm{K}_0(A)$. That is,
\[
([p]-n[1],\tau)\mapsto ([p]-n[1])^{\widehat{}}(\tau).
\]
Alternatively, in the stably projectionless case one may follow
the approach used in \cite{elliott3} and define the pairing 
of $[p]-n[1]$ and a lower semicontinuous trace $\tau$ by 
first finding $A_\tau\subseteq A$, a non-zero hereditary subalgebra
on which $\tau$ is bounded, subsequently finding a projection $p'\in \bigcup_n \mathrm M_n(A_\tau^\sim)$
Murray-von Neumann equivalent to $p$, and setting $([p]-n[1],\tau):=\tilde\tau(p')-n\tilde\tau(1)$,
where $\tilde\tau$ denotes the bounded extension of $\tau$ to $A_\tau^\sim$. 

\begin{corollary}\label{classsimple}
Let $A$ and $B$ be simple C$^*$-algebras that are expressible as inductive limits of 1-dimensional
NCCW complexes with trivial $\mathrm{K}_1$-group. Suppose that $\mathrm{K}_0(A)\cong \mathrm{K}_0(B)$ as ordered groups, $\T_0(A)\cong \T_0(B)$ as topological cones, and that these isomorphisms are compatible with the pairing between 
$\mathrm{K}_0$ and $\T_0$. Then $A\otimes\mathcal K$ and $B\otimes\mathcal K$ are isomorphic.

Furthermore, suppose that we have one of the following two cases:
\begin{enumerate}
\item[\textit{(a)}] $A$ and $B$ are both unital and the isomorphism from $\mathrm K_0(A)$
to $\mathrm K_0(B)$ maps $[1_A]$ to $[1_B]$, 

\item[\textit{(b)}] neither $A$ nor $B$ is unital and the isomorphism between $\T_0(A)$ and $\T_0(B)$
maps the tracial states of $A$ bijectively into the tracial states of $B$.
\end{enumerate}
Then $A$ and $B$ are isomorphic.
\end{corollary}
\begin{proof}
Notice that $A\otimes\mathcal K$ is projectionless if and only if $B\otimes\mathcal K$
is projectionless, as this property is equivalent to the  $\mathrm{K}_0$-groups having trivial cones
of positive elements (i.e., trivial order). Thus, either by Proposition \ref{Cusimprojections} or Proposition \ref{Cusimprojectionless},
depending on which is applicable, we have that $\Cu^\sim(A)\cong \Cu^\sim(B)$.
Moreover,  by  (a) or (b), again depending on the case at hand, the isomorphism from  $\Cu^\sim(A)$ to $\Cu^\sim(B)$ maps the class $[s_A]$ into the class $[s_B]$, where $s_A$
and $s_B$ are strictly positive elements of $A$ and $B$ respectively. It follows by
Corollary \ref{isomorphism} that $A$ and $B$ are isomorphic. 
\end{proof}

\begin{remark}\label{Cufails}
In Proposition \ref{Cuprojectionless} we have shown that if a C$^*$-algebra $A$ is stably projectionless and an inductive limit of 1-dimensional NCCW complexes  then $\Cu(A)$ is determined solely by $\T_0(A)$. It is known that among these C$^*$-algebras there are many that have non-trivial $\mathrm{K}_0$-group. In fact, it is explained  in \cite[Theorems 5.2.1 and 5.2.2]{elliott2} how such examples can be obtained as inductive limits of  1-dimensional NCCW complexes with trivial $\mathrm{K}_1$-group. It follows that while $\Cu^\sim$ is a classifying functor for these C$^*$-algebras, $\Cu$ is not, as it fails to account for their $\mathrm{K}_0$-groups.
\end{remark}

\subsection{Embeddings of $\mathcal Z$}
The Jiang-Su algebra is a simple unital exact non-elementary  C$^*$-algebra with a unique tracial state and strict comparison of positive elements (i.e., if $\widehat{[a]}\leq (1-\epsilon)\widehat{[b]}$ for some $\epsilon>0$, then $[a]\leq [b]$ for any Cuntz semigroup elements $[a]$ and $[b]$; see \cite{rordamZstable}).
These properties suffice to compute its Cuntz semigroup (see \cite{brown-toms}):
\[
\Cu(\mathcal Z)\cong \N\sqcup [0,\infty].
\]
More generally, if $A$ is a non-elementary simple unital C$^*$-algebra with a unique 2-quasitracial state and with strict comparison of positive elements, then 
\[
\Cu(A)\cong \mathrm V(A)\sqcup [0,\infty].
\]
Now assume that $A$ has stable rank one. We can then apply the classification result from the previous section with domain $\mathcal Z$
and with codomain the C$^*$-algebra $A$. Since $\mathcal Z$ is unital, we can use the functor $\Cu$, instead of $\Cu^\sim$, to  classify homomorphisms from $\mathcal Z$ (by Theorem \ref{permanence} (i)). Notice that there exists a unique  morphism in $\CCu$ from $\Cu(\mathcal Z)$ to $\Cu(A)$ such that $\N\ni 1\mapsto [1]\in \mathrm V(A)$ (since it must also map  $1\in [0,\infty]$ to $1\in [0,\infty]$). It follows that, up to approximate unitary equivalence, 
there exists a unique unital homomorphism from $\mathcal Z$ to $A$. This can be turned into a characterization of $\mathcal Z$.

\begin{proposition}\label{Zembeddings}
Consider the class $\mathcal C$ of C$^*$-algebras that are unital, simple, non-elementary, of stable rank one, with a unique 2-quasitracial state and with strict comparison of positive elements.
Then the Jiang-Su algebra is the unique C$^*$-algebra in $\mathcal C$ with the property that for any $A\in \mathcal C$ there exists, up to approximate unitary equivalence,   a unique unital embedding
from $\mathcal Z$ to $A$.
\end{proposition} 
\begin{proof}
We have already argued in the previous paragraph that $\mathcal Z$ embeds unitally and in a unique way -- up to approximate unitary equivalence -- into any C$^*$-algebra in $\mathcal C$.   Suppose that $\mathcal Z'$ is another C$^*$-algebra
in $\mathcal C$ with this  property. Then there exist unital embeddings  $\phi\colon \mathcal Z\to \mathcal Z'$ and $\psi\colon \mathcal Z'\to \mathcal Z$.
By the uniqueness of unital embeddings of $\mathcal Z$ in $\mathcal Z$, the endomorphism  $\psi\circ \phi$ is  approximately inner. Similarly, $\phi\circ\psi$ is also approximately inner. It follows by a standard intertwining argument that $\mathcal Z'\cong \mathcal Z$ (see, e.g., \cite{elliotttowards}).    
\end{proof}

Among the notable C$^*$-algebras in the class $\mathcal C$ of the previous proposition is $\mathrm C_r^*(\mathbb F_\infty)$, the reduced C$^*$-algebra of the free group with infinitely many generators. It is well known that $\mathrm C_r^*(\mathbb F_\infty)$ is simple and has a unique tracial state. 
Since $\mathrm C_r^*(\mathbb F_\infty)$ is exact, its 2-quasitraces are traces by Haagerup's Theorem. Thus,
$\mathrm C_r^*(\mathbb F_\infty)$ has a unique 2-quasi-tracial state. In \cite{dykema-haagerup-rordam},
Dykema, Haagerup, and R\o rdam show that $\mathrm C_r^*(\mathbb F_\infty)$ has stable rank one. Finally,
the following result due to R\o rdam (private communication) implies that $\mathrm C_r^*(\mathbb F_\infty)$ has strict comparison of positive elements and is thus in $\mathcal C$. 

\begin{proposition}[R\o rdam]\label{proprordam}
Let $(A_i,\tau_i)_{i=1}^\infty$ be an infinite sequence of unital C*-algebras $A_i$ and faithful tracial states $\tau_i$ on $A_i$. Suppose that for infinitely many indices $i\in \N$ the C*-algebra $A_i$ contains a unitary zero on $\tau_i$. Consider  the reduced free product \[(A,\tau):=\mathop{\ast}_{1\leq i<\infty} \,(A_i,\tau_i).\] 
If $a,b$ are positive elements in $A\otimes\mathcal K$ such that
$d_\tau(a)<d_\tau(b)$ then $a$ is Cuntz smaller than $b$. In particular, $A$ has strict comparison of positive elements.
\end{proposition}

\begin{proof}
The proof follows closely the proof of \cite[Theorem 2.1 (i)]{dykema-rordam},
where the same statement is made for projections instead of positive elements.
Assume without loss of generality that $a,b\in A\otimes \mathcal K$ are positive contractions.
Observe that $d_\tau(a)<d_\tau(b)$ implies that for every $\epsilon>0$ there is $\delta>0$
such that $d_\tau((a-\epsilon)_+)<d_\tau((b-\delta)_+))$. Thus, it suffices to show that 
$[(a-\epsilon)_+]\leq [(b-\delta)_+]$. For each $\epsilon>0$, $(a-\epsilon)_+$ is Cuntz equivalent to an element in $\mathrm M_n(A)$ for some $n\in \N$. This reduces the proof to the case that $a,b\in \mathrm M_n(A)$  for some $n\in \N$. Furthermore, viewing $A$ as the inductive limit
of the finite free products $\displaystyle{\mathop{\ast}_{1\leq i\leq m}} \,(A_i,\tau_i)$, with $m=1,2,\dots$, we can assume that $a$
and $b$ belong to the image in $\mathrm M_n(A)$ of the $n\times n$ matrices over one of such finite free products. Let $\epsilon>0$ and find $\delta>0$ such that $d_\tau((a-\epsilon)_+)<d_{\tau}((b-\delta)_+)$. 
We can argue as in \cite[Theorem 2.1 (i)]{dykema-rordam} that there exists a unitary $u\in \mathrm M_n(A)$ such that $u^*C^*(a)u$ and $C^*(b)$
are free with respect to the canonical extension of $\tau$ to a tracial state in $\mathrm M_n(A)$
(which we will continue to denote by $\tau$). Set $u^*(a-\epsilon)_+u=a'$ and $(b-\delta)_+=b'$.  Let $p_{a'},p_{b'}\in (\mathrm M_n(A))^{**}$  be the support projections of $a'$ and $b'$. From $d_{\tau}(a')<d_{\tau}(b')$ we get that $\tilde \tau(p_{a'})<\tilde \tau(p_{b'})$, where $\tilde \tau$ is a normal extension of $\tau$ to $(\mathrm M_n(A))^{**}$ (see \cite[Remark 5.3]{ortega-rordam-thiel}). Moreover, the fact that $a'$ and $b'$ are free implies that $p_{a'}$ and $p_{b'}$ are free with respect to $\tilde\tau$. (To get this, use that $(a')^{1/n}$ and $(b')^{1/n}$ are free and pass to the limit in the formula expressing their mixed moments in terms of the moments of $(a')^{1/n}$ and of $(b')^{1/n}$.) By 
\cite[Proposition 1.1]{dykema-rordam1}, we have that $\|p_{a'}(1-p_{b'})\|=\alpha<1$.
Now from $(1-\alpha^2)p_{a'}\leq p_{a'}p_{b'}p_{a'}$ we get that \[
(1-\alpha^2)a'\leq (a')^{1/2} p_{b'}(a')^{1/2}\leq (a')^{1/2}f(b)(a')^{1/2}.\] 
Here $f(b)\in C^*(b)^+$ is chosen such that $f(b)b'=b'$. We conclude that
$[(a-\epsilon)_+]=[a']\leq [b]$ for every $\epsilon>0$. Thus, $[a]\leq [b]$.
\end{proof}

It is not known whether the C$^*$-algebras $\mathrm C_r^*(\mathbb F_n)$ have strict comparison of positive elements for $n<\infty$. Notice, however,  that since $\mathrm C_r^*(\mathbb F_\infty)$ embeds unitally in $\mathrm C_r^*(\mathbb F_n)$ for all $n\geq 2$, we have that $\mathcal Z$ embeds unitally in $\mathrm C_r^*(\mathbb F_n)$ for $n=2,3,\dots,\infty$.

\section{The crossed products $\mathcal O_2\rtimes_\lambda \R$}\label{crossed}
Let $\lambda\in \R$. Consider the action of $\R$ on the Cuntz algebra $\mathcal O_2$ given by
\[
\sigma_\lambda^t(v_1):=e^{2\pi i t}v_1,\, \sigma_\lambda^t(v_2):=e^{2\pi i \lambda t}v_2,
\]
where $v_1$ and $v_2$ are the partial isometries generating $\mathcal O_2$. These actions, and the resulting
crossed-product C$^*$-algebras $\mathcal O_2\rtimes_\lambda \R$, were first studied by Evans and subsequently by
Kisihimoto, Kumjian, Dean, and several other authors (see \cite{evans}, \cite{kishimoto-kumjian}, \cite{dean}). In \cite{kishimoto}, 
Kishimoto shows that if $\lambda<0$ 
the crossed product is purely infinite, while if $\lambda>0$ it is stably finite. In \cite{dean}, Dean computes
$\mathcal O_2\rtimes_\lambda \R$ for $\lambda\in \Q^+$ and uses this to conclude that for $\lambda>0$ in a dense $G_\delta$ subset of $\R^+$ 
that contains  $\Q^+$ the C$^*$-algebras $\mathcal O_2\rtimes_\lambda \R$ are inductive limits
of 1-dimensional NCCW complexes. Here we will show that the 1-dimensional NCCW complexes obtained
by Dean all have trivial $\mathrm K_1$-group. Since for $\lambda>0$ irrational the C$^*$-algebras 
$\mathcal O_2\rtimes_\lambda \R$ are simple, stable, projectionless, and with a unique trace,  we shall conclude by Corollary
\ref{classsimple} that these C$^*$-algebras are all isomorphic for 
a dense set of irrational numbers $\lambda>0$.

Let $p,q>0$ be relatively prime natural numbers. In \cite{dean}, Dean shows that there is a simple, stable,  
AF algebra $A(p,q)$, and an automorphism $\alpha\colon A(p,q)\to A(p,q)$ such that
$\mathcal O_2\rtimes_{p/q} \R$ is isomorphic to the mapping torus $M_{\alpha,p,q}$ of $(A(p,q),\alpha)$. That is,
\[
\mathcal O_2\rtimes_{p/q}\R\cong M_{\alpha,p,q}:=\{f\in \mathrm C([0,1],A(p,q))\mid f(0)=\alpha(f(1))\}.
\]
A description in terms of generators and relations of the algebra $A(p,q)$, and of the action of the automorphism $\alpha$, 
is given in \cite[Theorem 3.1]{dean}. The computation of the
K-theory of $A(p,q)$ is sketched at the end of \cite{dean} (and carried out in detail for $p=1$, $q=2$). 
From the discussion given in the final remarks of \cite[Section 5]{dean}, one can gather the following:
\begin{enumerate}
\item
There is an increasing sequence of finite dimensional algebras $(D_n)_{n=1}^\infty$
such that $A(p,q)=\overline{\bigcup_n D_n}$.

\item
For each $n$, $\mathrm K_0(D_n)\cong \Z^q$ and the inclusion $D_n\subseteq D_{n+1}$ induces
in $\mathrm K_0$ a map $\Z^q\to \Z^q$ given by the matrix $A^{q+1}$, with
\[
A=
\begin{pmatrix}
 & 1 &&&\\
&& 1 &&\\
&&&\ddots&\\
&&&&1\\
&&1&\dots&1
\end{pmatrix}.
\]
The last row of $A$ has 1 in the entries $A_{p-q,q}$ and $A_{q,q}$ and zeroes elsewhere.

\item
For every $n$, $\alpha(D_n)\subseteq D_{n+1}$ and the homomorphism $\alpha|_{D_n}\colon D_n\to D_{n+1}$ induces in
$\mathrm K_0$ the map $A^{q}\colon \Z^q\to A^q$, with $A$ as above.
\end{enumerate}

The above points may be summarized by the commutative diagram
\[
\xymatrix{
D_1\ar@{^{(}->}[r]\ar@{^{(}->}^{\alpha|_{D_1}}[d] &
D_2\ar@{^{(}->}[r] \ar@{^{(}->}^{\alpha|_{D_2}}[d]& \dots & A(p,q)\ar[d]^\alpha\\
D_2\ar@{^{(}->}[r] & D_3\ar@{^{(}->}[r] & \dots & A(p,q)
},
\]
and the diagram induced in $\mathrm K_0$,
\[
\xymatrix{
\Z^q\ar[r]^{A^{q+1}}\ar[d]^{A^{q}} & \Z^q\ar[r]^{A^{q+1}} \ar[d]^{A^{q}}& \dots & \Z^q\ar[d]^{A^q}\\
\Z^q\ar[r]^{A^{q+1}} & \Z^q\ar[r]^{A^{q+1}} & \dots & \Z^q.
}
\]

The mapping torus $M_{\alpha,p,q}$ is the inductive limit of the 1-dimensional NCCW complexes
\[
M_{n,p,q}=\{f\in \mathrm C([0,1],D_{n+1})\mid f(1)\in D_n,\, f(0)=\alpha(f(1))\}.
\]
By Lemma \ref{K10}, $\mathrm K_1(M_{n,p,q})=0$ if and only if $\mathrm K_0(\phi_0)-\mathrm K_0(\phi_1)$
is surjective. In this case, $\phi_0$ is the inclusion of $D_n$ in $D_{n+1}$, while $\phi_1$ is the
restriction of $\alpha$ to $D_n$ (with codomain $D_{n+1}$). Identifying $\mathrm K_0(D_n)$ and $\mathrm K_0(D_{n+1})$ with $\Z^q$, we have
$\mathrm K_0(\phi_0)=A^{q+1}$ and $\mathrm K_0(\phi_1)=A^q$.
Thus,
\[
\mathrm K_0(\phi_0)-\mathrm K_0(\phi_1)=A^{q+1}-A^q=A^q(A-I).
\]
Since the characteristic polynomial of $A$ is $t^q-t^{q-p}-1$, we see that $\det(I-A)=-1$ and $\det(-A)=-1$. Thus, $A^{q+1}-A^q$ is an invertible map from $\Z^q$ to $\Z^q$, and so $\mathrm{K}_1(M_{n,p,q})=0$.
By Theorem \ref{main}, the functor $\Cu^\sim$ classifies homomorphisms from  the NCCW 
complexes  $M_{n,p,q}$, and from their sequential inductive limits. By \cite[Theorem 4.10]{dean} (or rather, by its proof),
these inductive limits include all C$^*$-algebras $\mathcal O_2\rtimes_\lambda \R$ for a dense $G_\delta$ set of positive numbers $\lambda$ 
that includes $\Q^+$.

\begin{corollary}
The crossed products
$\mathcal O_2\rtimes_\lambda \R$ are all isomorphic
for $\lambda$ belonging to a dense subset of $\R^+\backslash \Q$ of second Baire category.
\end{corollary}
\begin{proof}
By the discussion above,
for a $G_\delta$ set  $\Lambda\subseteq \R^+$  the crossed products 
$\mathcal O_2\rtimes_\lambda \R$
are inductive limits of 1-dimensional  NCCW complexes with trivial $\mathrm K_1$-group.
On the other hand, for $\lambda>0$ irrational, these crossed products are simple, stable, 
have trivial $\mathrm K_0$-group, and tracial cone $\R^+$ (see \cite{kishimoto-kumjian}). It follows by Corollary \ref{classsimple}
that for $\lambda\in \Lambda\backslash \Q$  the crossed products $\mathcal O_2\rtimes_\lambda \R$ are all isomorphic.
\end{proof}

\begin{remark}
Quasifree actions may also be defined on $\mathcal O_n$ and they have been studied by the authors
cited above. However, for the resulting crossed products the $\mathrm K_1$-group is $\Z/(n-1)\Z$. 
So their classification lies beyond the scope of the results obtained here.
\end{remark}

\end{document}